\documentclass{article}
\usepackage{a4,amssymb,amsmath,color,graphicx}
\usepackage{babel}

\def\eins{\mbox{1\hskip-0.24em l}}
\def\T{^{\sf T}}
\def\mT{^{\sf -T}}

\newcommand{\N}{ {\mathbb N} }
\newcommand{\R}{ {\mathbb R} }
\newcommand{\MM}{{\mathbb M}}

\newcommand{\CC}{{\cal C}}

\newcommand{\LL}{{\cal L}}
\newcommand{\PP}{{\cal P}}
\newcommand{\QQ}{{\cal Q}}
\newcommand{\Sc}{{\cal S}}
\newcommand{\diag}{\,\mbox{diag}}

\newcommand{\cc}{{\bf c}}
\newcommand{\pts}{:} 
\newcommand{\gdw}{\ \iff\ }
\newcommand{\folgt}{\ \Rightarrow\ }

\newcommand{\qed}{\qquad\mbox{$\square$}}
\newcommand{\mynegspace}{\hspace{-0.12em}}
\newcommand{\ltnorm}{\rvert\mynegspace\rvert\mynegspace\rvert}
\newcommand{\rtnorm}{\rvert\mynegspace\rvert\mynegspace\rvert}

\newtheorem{remark}{Remark}[section]
\newtheorem{theorem}{Theorem}[section]
\newtheorem{lemma}{Lemma}[section]

\newtheorem{corol}{Corollary}[section]

\begin{document}
\title{Variable-Stepsize Implicit Peer Triplets\\
in ODE Constrained Optimal Control}
\author{Jens Lang \\
{\small \it Technical University Darmstadt,
Department of Mathematics} \\
{\small \it Dolivostra{\ss}e 15, 64293 Darmstadt, Germany}\\
{\small lang@mathematik.tu-darmstadt.de} \\ \\
Bernhard A. Schmitt \\
{\small \it Philipps-Universit\"at Marburg,
Department of Mathematics,}\\
{\small \it Hans-Meerwein-Stra{\ss}e 6, 35043 Marburg, Germany} \\
{\small schmitt@mathematik.uni-marburg.de}}
\maketitle

\begin{abstract}
This paper is concerned with the theory, construction and application of
implicit Peer two-step methods that are super-convergent for variable stepsizes, i.e., preserve their
classical order achieved for uniform stepsizes when applied to ODE constrained optimal control problems
in a first-discretize-then-optimize setting.
We upgrade our former implicit two-step Peer triplets constructed in [Algorithms, 15:310, 2022] to get
ready for dynamical systems with varying time scales without loosing efficiency.
Peer triplets consist of a standard Peer method for interior time steps supplemented by matching methods 
for the starting and end steps.
A decisive advantage of Peer methods is their absence of order reduction since they use stages of the same high stage order.
The consistency analysis of variable-stepsize implicit Peer methods results in additional order conditions and severe new difficulties for uniform zero-stability, which intensifies the demands on the Peer triplet.
Further, we discuss the construction of 4-stage methods with order pairs (4,3)
and (3,3) for state and adjoint variables in detail and provide four Peer triplets of practical interest. We
rigorously prove convergence of order $s-1$ for $s$-stage Peer methods applied on grids
with bounded or smoothly changing stepsize ratios.
Numerical tests show the expected order of convergence for the new variable-stepsize Peer triplets.
\end{abstract}


\section{Introduction}
Recently, we have developed and tested third- and fourth-order implicit Peer two-step methods \cite{LangSchmitt2022a,LangSchmitt2022b,LangSchmitt2023a}
to solve ODE constrained optimal control problems of the form
\begin{align}
\mbox{minimize } C\big(y(T)\big) \label{OCprob_objfunc} &\\
\mbox{subject to } y'(t) =& \,f\big(y(t),u(t)\big),\quad
u(t)\in U_{ad},\;t\in(0,T], \label{OCprob_ODE}\\
y(0) =& \,y_0, \label{OCprob_ODEinit}
\end{align}
with the state $y(t)\in\R^m$, the control $u(t)\in\R^d$,
$f: \R^m\times\R^d\mapsto\R^m$, the objective function $C: \R^m\mapsto\R$, where the set of 
admissible controls $U_{ad}\subset\R^d$ is closed and convex. The design of efficient time integrators for the numerical solution of 
such problems with large $m$ arising from semi-discretized time-dependent partial differential equations
is still of great interest since difficulties arise through additional adjoint order conditions, some recent literature is \cite{AlbiHertyPareschi2019,AlmuslimaniVilmart2021}.
Implicit Peer two-step methods overcome the structural disadvantages of one-step and multi-step
methods such as symplectic or generalized partitioned Runge-Kutta methods \cite{BonnansLaurentVarin2006,HertyPareschiSteffensen2013,LiuFrank2021} and backward differentiation formulas \cite{BeigelMommerWirschingBock2014}.
They avoid order reduction of one-step methods, e.g., for boundary control problems of PDEs \cite{LangSchmitt2023a}, and have good stability properties.
Peer methods also allow the approximation of adjoints in a first-discretize-then-optimize (FDTO) approach with higher order
which still seems to be an unsolved problem for multi-step methods.
FDTO is the most commonly used method and possesses the advantage of providing consistent gradients for state-of-the-art optimization algorithms. We refer the reader to the detailed discussions in our previous papers
\cite{LangSchmitt2022a,LangSchmitt2022b,LangSchmitt2023a}. 
\par
So far we have considered the use of constant stepsizes to approximate the ODE system \eqref{OCprob_ODE}. However,
dynamical systems with sub-processes evolving on many different time scales are
ubiquitous in applications. Their efficient solution is greatly enhanced by automatic time step variation. 
A popular strategy is to focus on adaptively optimizing time grids in accordance with local error control, i.e.,
errors made within a single integration step. However, in optimal control algorithm, the approximation property of
adjoint variables and controls heavily depend on the global state errors at any 
discrete time point. Ideally, higher-order time integrators should not suffer from order reduction when they are applied
with smoothly varying stepsizes. 
Since the optimal control problem requires the solution of a boundary value problem, where solutions from the whole grid have to be saved, higher order methods and adaptive grids may lead to considerable savings in memory.
So this paper is concerned with the theory, construction and application of implicit Peer two-step methods that are 
super-convergent for variable stepsizes, i.e., preserve their classical order achieved for uniform stepsizes.
In fact, we are considering triplets of methods where some common standard Peer method is accompanied by different start and end methods satisfying appropriate matching conditions.
\par
Introducing for any $u\in U_{ad}$ the normal cone mapping
\begin{align}
\label{def_cone}
N_U(u) =&\, \{ w\in\R^d: w^T(v-u)\le 0 \mbox{ for all } v\in U_{ad}\},
\end{align}
the first-order Karush–Kuhn–Tucker (KKT) optimality conditions read \cite{Hager2000,Troutman1996}
\begin{align}
y'(t) =& \,f\big(y(t),u(t)\big),\quad t\in(0,T],\quad y(0)=y_0, \label{KKT_state}\\
p'(t) =& \,-\nabla_y f\big(y(t),u(t)\big)\T p(t),\quad t\in[0,T),
\quad p(T)=\nabla_y {\CC}\big(y(T)\big)\T, \label{KKT_adj}\\
& \,-\nabla_u f\big(y(t),u(t)\big)\T p(t) \in N_U\big(u(t)\big),\quad t\in[0,T]. \label{KKT_ctr}
\end{align}
Throughout the paper, we will make the following two assumptions: (i) There exists a local solution 
$(y^\star,u^\star,p^\star)$ of the KKT system and (ii) the control uniqueness property as stated
in \cite{Hager2000} is valid, i.e.,
for any $(y,p)$ sufficiently close to $(y^\star,p^\star)$, there exists a locally unique minimizer 
$u=u(y,p)$ of the Hamiltonian $H(y,u,p):=p\T f(y,u)$ over all $u\in U_{ad}$. Then 
the control $u(t)$ can be explicitly eliminated, yielding the reduced boundary value problem
\begin{align}
\label{RWPy}
y'(t)=&\,g\big(y(t),p(t)\big),\quad y(0)=y_0,\\
\label{RWPp}
p'(t)=&\,\phi\big(y(t),p(t)\big),\quad p(T)=\nabla_y C\big(y(T)\big)\T,
\end{align}
with the source functions defined by
\begin{align}
g(y,p) :=& f\big(y,u(y,p)\big),
\quad \phi(y,p) := -\nabla_yf\big(y,u(y,p)\big)\T p.
\end{align}
This boundary value problem is used in our consistency and convergence analysis.
\par
The paper is organised as follows. In Chapter~\ref{sec:ipeer}, we derive the discrete
system of Karush-Kuhn-Tucker conditions for Peer triplets and reformulate it as
boundary value problem with eliminated control. A consistency analysis for the standard and
the boundary methods is presented in Chapter~\ref{sec:conistency}. Uniform zero-stability 
is established by enforcing boundedness of matrix families in Chapter~\ref{SBMF}. 
In Chapter~\ref{SDesign4}, we will design new 4-stage Peer triplets of practical interest.
Convergence for variable stepsizes is studied in Chapter~\ref{sec:convergence}. We discuss 
the results of two optimal control problems in Chapter~\ref{sec:tests} and conclude with a summary
in Chapter~\ref{sec:summary}.

\section{Implicit Peer two-step methods with variable stepsizes}\label{sec:ipeer}
We follow the \textit{first-discretize-and-then-optimize} approach.
In every time step from $t_n$ to $t_{n+1}:=t_n+h_n$, two-step Peer methods use $s$ stages $Y_{ni}\cong y(t_{ni}),\,i=1,\ldots,s$, 
as solution approximations at off-step points $t_{ni}=t_n+h_nc_i,\,i=1,\ldots,s$, based on a fixed set of nodes which 
may be considered as elements of a node vector $\cc=(c_1,\ldots,c_s)$.
For ease of writing, also the stages are collected in a large vector, e.g.,  $Y_n=\big(Y_{ni}\big)_{i=1}^s$.
For variable stepsizes, the two-step structure of the Peer methods requires that some of its coefficients have to change between steps.
Hence, we associate with every time step $t_n\to t_{n+1}$ of the Peer methods given in a redundant formulation by
\begin{align}\label{PeerStd}
A_nY_n=B_nY_{n-1}+h_n K_nF(Y_n,U_n),\ n=1,\ldots,N,
\end{align}
for the approximation of \eqref{OCprob_ODE}, an individual set of three coefficient matrices $(A_n,B_n,K_n)$.
Here, $F(Y_n,U_n)=\big(f(Y_{ni},U_{ni})\big)_{i=1}^s$.
As in \eqref{PeerStd}, we will use for a simple matrix like $A_n\in\R^{s\times s}$ and for its 
Kronecker product $A_n\otimes I_m$ with the identity matrix the same symbol.
For the stage approximations $Y_0$ of the first interval, a one-step method of Runge-Kutta type
\begin{align}\label{PeerStrt}
A_0Y_0=a\otimes y_0+h_0 K_0F(Y_0,U_0),
\end{align}
with $a\in\R^s$ is considered. We allow for a general linear combination for the solution at the end point
\begin{align}\label{PeerEnd}
 y_h(T):=\sum_{j=1}^s w_j Y_{nj}=\left(w\T\otimes I_m\right) Y_N,
\end{align}
in order to avoid the necessity of choosing $c_s=1$. The matrices $A_n$ are assumed to be nonsingular.
For reasons of computational efficiency, lower triangular matrices $A_n$ and diagonal matrices $K_n$ 
will be used with the exception of the boundary steps $n=0,N$.
\par
Considering the discrete Lagrange function
\begin{align}
\notag
L(Y,U,P)
&=\,C(y_h(T))\,-P_0\T\big(A_0 Y_0-a\otimes y_0-h_0 K_0F(Y_0,U_0)\big)\\
\label{LgrFkG}
&\,-\sum_{n=1}^{N} P_n\T\big(A_nY_n-B_nY_{n-1}-h_nK_nF(Y_n,U_n)\big)
\end{align}
with multipliers $P_n\in\R^{ms}$, the additional discrete KKT conditions read, see \cite{LangSchmitt2022b},
\begin{align}\label{AdjABK}
A_n\T P_n=&\,B_{n+1}\T P_{n+1}+h_n\nabla_YF(Y_n,U_n)\T K_n\T P_n,\ 0\le n\le N-1,\\
\label{AdjStrt}
A_N\T P_N=&\,w\otimes p_h(T)+h_N\nabla_YF(Y_N,U_N)\T K_N\T P_N,\ n=N,
\\\label{KKTU}
&-\nabla_U F(Y_n,U_n)\T K_n\T P_n\in N_{U^*}(U_n),\ 0\le n\le N,
\end{align}
where $p_h(T)=\nabla_y C(y_h(T))\T$ and 
\begin{align*} 
N_{U^*}(u)=\left\{ w\in\R^{sd}: w^T(v-u)\le 0 \mbox{ for all } v\in U_{ad}^s\subset\R^{sd}\right\}.
\end{align*}
Hence, the KKT conditions constitute a discrete boundary value problem for $(Y,U,P)$ described by the 5 equations \eqref{PeerStd}, \eqref{PeerStrt}, \eqref{AdjABK}, \eqref{AdjStrt}, \eqref{KKTU}.
We note that the forward step \eqref{PeerStd} and the adjoint step \eqref{AdjABK} may be simplified, e.g., by multiplication with $A_n^{-1}$ or $A_n\mT$, respectively.
However, these changes lead to different \glqq simplified\grqq{} coefficients.
Hence, the redundant formulation of the standard method with 3 matrices $(A_n,B_n,K_k)$ adds degrees of freedom required for the additional order conditions, \cite{LangSchmitt2022a}.
\par
Eliminating the discrete controls $U_n=U_n(Y_n,P_n)$ by the control uniqueness property and defining
\begin{align}
\Phi(Y_n,K_n\T P_n):=&\,\left( \phi(Y_{ni},(K_n\T P_{n})_i)\right)_{i=1}^{s},
\quad G(Y_n,P_n):=\left( g(Y_{ni},P_{ni})\right)_{i=1}^{s},
\end{align}
the discrete boundary value problem reads in compact form,
\begin{align}
A_0Y_0=&\,a\otimes y_0+h_0K_0G(Y_0,P_0),\label{dBVPY0}\\[1mm]
A_nY_n=&\,B_nY_{n-1}+h_nK_nG(Y_n,P_n),\quad 1\le n\le N,\label{dBVPYn}\\[1mm]
y_h(T)=&\,(w\T\otimes I_m)Y_N,\label{dBVPyT}\\[1mm]
A_n\T P_n=&\,B_{n+1}\T P_{n+1}-h_n\Phi(Y_n,K_n\T P_n),\ 0\le n\le N-1,\label{dBVPPn}\\[1mm]
A_N\T P_N=&\,w\otimes \nabla_y C(y_h(T))\T-h_N\Phi(Y_N,K_N\T P_N),\ n=N.\label{dBVPPN}
\end{align}
We will now study the consistency of the overall scheme.

\section{Consistency}\label{sec:conistency}
\subsection{Order conditions for the inner grid}
Since we are interested in localized error estimates taking account of nodes outside the standard interval $[0,1]$, i.e. with $c^\ast:=\max\{c_1,\ldots,c_s,1\}\ge1$, we consider the whole interval $[0,T]$ as part of the union $\bigcup_{n=0}^N\pi_n\subset[0,T^\ast]$ of subintervals $\pi_n:=[t_n,t_n+h_nc^\ast]$.
The order conditions from earlier papers \cite{LangSchmitt2022a,LangSchmitt2022b} have to be extended now with the aid of the scaling matrix $S_{n,r}=\diag_i(\sigma_n^{i-1})\in\R^{r\times r}$, $1\le r\le s$, depending on the local stepsize ratio $\sigma_n=h_n/h_{n-1}>0$, $n\ge 1$, which will be restricted by lower and upper bounds later on.
In principle, we will again consider conditions with two orders $r,q\le s$ applied uniformly for all $\sigma\in\R$, where $r$ denotes the local order of the standard scheme \eqref{PeerStd} and $q$ the local order for the adjoint scheme \eqref{AdjABK}.
For the forward step \eqref{PeerStd}, such conditions have already been presented in, e.g. \cite{SchmittWeinerErdmann2005}.
They will be given in \eqref{OBvStd}.
Hence we derive in detail only the order conditions for the adjoint step \eqref{AdjABK} with a diagonal matrix $K_n$.
Using Taylor expansion at $t_n$ for $t_{nj}=t_n+h_nc_j$ and $t_{n+1}+h_{n+1}c_j=t_n+h_n(1+\sigma_{n+1}c_j)$, this step has the following residuals $\tau_n^P$ if it is used with values $y(t_{ni}),p(t_{ni})$ and $\nabla_y f(y(t_{ni}),u(t_{ni}))\T p(t_{ni})=-p'(t_{ni})$.
Let $A_n:=(a_{ij})_{i,j=1}^s$, $B_n:=(b_{ij})_{i,j=1}^s$ and $K_n=\diag(\kappa_{jj})_{j=1}^s$. With $z_n:=h_n d/dt$ and $\exp_q(z)=\sum_{j=0}^{q-1}z^j/j!$, we then obtain for the local error the expression
\begin{align}\notag
 A_n\T\tau^P_{n}:=&\left(\sum_{i=1}^s\big(a_{ij}p(t_{ni})-b_{ij}p(t_{n+1,i})\big)+
 h_n\kappa_{jj}p'(t_{nj})\right)_{j=1}^s
 \\\label{lokFep}
 =&\Big(A_n\T\exp_q(\cc z_n)-B_n\T\exp_q\big((\eins+\sigma_{n+1}\cc)z_n\big)+z_nK_n\exp_{q-1}(\cc z_n)\Big)p|_{t_n}
 +O(z_n^q p).
\end{align}
We introduce the powers of the node vector $\cc^k=(c_1^k,\ldots,c_s^k)\T$, the Vandermonde matrices $V_q=(\eins,\cc,\ldots,\cc^{q-1})\in\R^{s\times q}$, the Pascal matrix $\PP_q=\big({j-1\choose i-1}\big)\in\R^{q\times q}$ and the nilpotent matrix $\tilde E_q=(i\delta_{i+1,j})\in\R^{q\times q}$ commuting with $\PP_q=\exp(\tilde E_q)$.
Then, with the identity
\[ \big((1+\sigma c_i)^{k-1}\big)_{i,k=1}^{s,q}=\sum_{l=1}^kc_i^{l-1}\sigma^{l-1}{k-1\choose l-1}=V_q S_q(\sigma)\PP_q,
\]
the conditions for adjoint local order $q\le s$ may be deduced from \eqref{lokFep} and combined in matrix form to
\begin{align}\label{OBast}
 A_n\T V_q=&B_{n+1}\T V_qS_{n+1,q}\PP_q-K_nV_q\tilde E_q.
\end{align}
Transposing equation \eqref{OBast}, we may present the conditions for local order $r\le s$ of the forward method and local order $q\le s$ of the adjoint method together as
\begin{align}\label{OBvStd}
 A_nV_r=&B_nV_r\PP_r^{-1}S_{n,r}^{-1}+K_nV_r\tilde E_r,
 \\\label{OBaStd}
 V_q\T A_n=&\PP_q\T S_{n+1,q}V_q\T B_{n+1}-\tilde E_q\T V_q\T K_n.
\end{align}
The equations of lowest order, often called preconsistency conditions will be used frequently and are repeated here for easier reference,
\begin{align}\label{precons}
 A_n\eins=B_n\eins,\quad \eins\T A_n=\eins\T B_{n+1}.
\end{align}
Also for later reference, we note the precise form of the local error resulting from these order conditions.
Since the adjoint step analyzed in \eqref{lokFep} uses nodes in $\pi_n\cup\pi_{n+1}$, $n<N$, and the forward step from $\pi_{n-1}\cup\pi_n$, we introduce the following short-hand notation for the norms used in the remainder terms.
For $v\in C[0,T^\ast]$ let
\begin{align}
 \|v\|_{[n]}:=\max\{\|v(x)\|:\ x\in\pi_{n-1}\cup\pi_n\cup\pi_{n+1}\},
 \ 0\le n\le N,
\end{align}
with the agreement that $\pi_0=\pi_{N+1}=\emptyset$.
Now, if \eqref{OBvStd}, \eqref{OBaStd} hold and $y,p\in C^{k+1}[0,T^\ast]$ for $k\le r$ resp. $k\le q$, then the local errors satisfy
\begin{align}\label{lokFRestv}
\tau_n^Y=&\;h_n^k\beta_k(\sigma_n)y^{(k)}(t_n)+O\big(h_n^{k+1}\|y^{(k+1)}\|_{[n]}\big),\\\notag
 &\beta_k(\sigma):=
\frac1{k!}A_n^{-1}\big(A_n\cc^k-B_n(\cc-\eins)^k\sigma^{-k}-kK_n\cc^{k-1}\big),
\\\label{lokFResta}
\tau^P_{n}=&\;h_n^k\beta_k^\dagger(\sigma_{n+1})p^{(k)}(t_n)
 +O\big(h_n^{k+1}\|p^{(k+1)}\|_{[n]}\big),\\
 &\beta_k^\dagger(\sigma):=\frac1{k!} A_n\mT\big(A_n\T\cc^k -B_{n+1}\T(\eins+\sigma\cc)^k+kK_n\cc^{k-1}\big).
\end{align}
Using the current stepsize $h_n$ also in the remainder is justified since stepsizes will be restricted to some interval $\sigma_n\in[\underline\sigma,\bar\sigma]$ on the positive axis, for instance with $\underline\sigma=1/\bar\sigma<1$.
Of course, only the remainder terms appear here for  $k<r$ resp. $k<q$.
\par
Before proceeding, we like to recall our basic design strategy from \cite{LangSchmitt2022b} for matching order conditions for the three parts of the Peer triplets.
Since the boundary steps \eqref{PeerStrt} and \eqref{AdjStrt} are applied once only, they contribute to the global error with their local order.
For the time steps of the inner grid, we will prevent the usual loss of one order from the local to the global error by exploiting again a super-convergence property.
Accordingly, all steps forward in time have to satisfy conditions with the same local order $r<s$  and the adjoint steps for order $q<s$ for arbitrary stepsize ratios $\sigma_n$.
\par
Obviously, the two boundary steps \eqref{PeerStrt}, \eqref{AdjStrt} have one-step character and a triangular form of, e.g., $A_0,K_0$ would limit the (local) stage order to two.
Hence, we have to consider coefficient matrices $A_n,K_n$, $n=0,N,$ of more general form for the two boundary methods.
Of course, the loss of triangularity increases the computational expense by much.
However, this affects only the boundary steps which present only a very small part of the full boundary value problem.
In the boundary steps, $n=0,N,$ one additional {\em one-leg-condition}
\begin{align}\label{OBU}
(\cc^{q_b-1})\T K_n=\eins\T K_n C^{q_b-1},\quad C:=\diag(c_1,\ldots,c_s),
\end{align}
applies with order $q_b=q-1$, which is nontrivial for $q=3$, see (61) in \cite{LangSchmitt2022b}.
It is obviously satisfied with any $q_b\in\N$ for diagonal matrices $K_n$ like those from the rest of the grid.

\subsection{Step-independent matrices $A,K$ in the standard method}\label{SAKconst}
Using different coefficients in every time step would mean a tremendous overhead for the imple\-mentation of such methods.
Instead, we will use exceptional coefficients in the boundary steps only and use fixed coefficients $A_n\equiv A$, $K_n\equiv K$ of our {\em standard method} in the inner time steps with $1\le n<N$.
There, only $B_n$ may change with $n$.
In fact, the order conditions will show that it will depend on the stepsize ratio $\sigma_n$ only: $B_n=B(\sigma_n),\,n=1,\ldots,N$, with a certain matrix function $B(\sigma)$ depending on the parameter $\sigma\in\R$.
This strategy has been used earlier, e.g.  \cite{SchmittWeinerErdmann2005}, but now it has far reaching  consequences on the structure of all three coefficient matrices $A,B(\sigma),K$.
\par
Since now 2 of the 3 matrices in the order conditions \eqref{OBvStd} are independent of $\sigma_n$, the same holds for the third one containing $B(\sigma_n)$, which means that
\begin{align}\label{OBvsw}
 AV_r-KV_r\tilde E_r =B(\sigma_n)V_r\PP_r^{-1}S_{n,r}^{-1}
 \equiv B(1)V_r\PP_r^{-1},
\end{align}
since $S_{n,r}=I_r$ for $\sigma_n=1$.
Hence, $A$ and $K$ may be dropped in the present discussion.
The same argument for condition \eqref{OBaStd} means that also
\begin{align}\label{OBasw}
\PP_q\mT (V_q\T A+\tilde E_q\T V_q\T K)=S_{n,q}V_q\T B(\sigma_{n})
 \equiv V_q\T B(1)
\end{align}
is independent of $\sigma_n$.
Combining both conditions fixes large parts of the matrix function $B(\sigma)$ and as a consequence also of $A$ and $K$.
We note that the matrix $\QQ_{q,r}$ appearing in the next lemma plays a fundamental role in the matching conditions for the three methods constituting the triplet, see also \cite{LangSchmitt2022b}.
\par
It is obvious that due to the redundant formulation of the standard Peer method with three matrices $(A,B_n,K)$ any nontrivial multiple of it also satisfies the order conditions.
In order to eliminate unnecessary degrees of freedom from the upcoming discussion which will disappear later on we assume the following normalization:
\begin{align}\label{eAe}
 \eins\T A\eins =1,
\end{align}
which will naturally come up through the order conditions for the boundary steps.
\begin{lemma}\label{LBsig}
Let the order conditions \eqref{OBvStd} and \eqref{OBaStd} be satisfied with $1\le q,r\le s$ for two different stepsize ratios $\sigma_n$, at least, and assume \eqref{eAe}.
Then the matrix $\QQ_{q,r}(\sigma):=V_q\T B(\sigma) V_r\PP_r^{-1}$ satisfies
\begin{align}\label{QQform}
 \QQ_{q,r}(\sigma)
 =e_1e_1\T\in\R^{q\times r}
\end{align}
and is independent of $\sigma$.
\end{lemma}
{\bf Proof:}
By \eqref{precons} the assumption \eqref{eAe} leads to $1=\eins\T A\eins=\eins\T B(\sigma)\eins=e_1\T\QQ_{q,r}(\sigma)e_1$.
Multiplying \eqref{OBvsw} by $V_q\T$ from the left and \eqref{OBasw} by $V_r\PP_r^{-1}$ from the right yields
\begin{align*}
V_q\T B(\sigma_n)V_r\PP_r^{-1}S_{n,r}^{-1}=\QQ_{q,r}(1)\gdw \QQ_{q,r}(\sigma_n)=\QQ_{q,r}(1)S_{n,r},\\
S_{n,q}\QQ_{q,r}(\sigma_n)=\QQ_{q,r}(1)\folgt
 \QQ_{q,r}(1) = S_{n,q}\QQ_{q,r}(\sigma_n)=S_{n,q}\QQ_{q,r}(1)S_{n,r},
\end{align*}
where the identity from the first line was inserted in the last step.
Hence the elements of $\QQ_{q,r}(1)=(\gamma_{ij})$ satisfy $\gamma_{ij}(1-\sigma_n^{i+j-2})=0$ for $i\le q,\,j\le r$.
Fulfilling this condition for two different values of $\sigma_n$ leaves as only solution the one from assertion \eqref{QQform} since $\gamma_{11}=1$.
\qed
\par
This Lemma is the first place showing that matrices congruent to the original coefficients $A,B(\sigma),K$ may have a very restricted form.
Since these matrices will be discussed frequently, we introduce a shorthand notation for them and their submatrices:
\begin{align}\label{DachMat}
 \hat A_{q,r}:=V_q\T A V_r,\quad
 \hat B_{q,r}(\sigma):=V_q\T B(\sigma) V_r,\quad
 \hat K_{q,r}:=V_q\T K V_r,
\end{align}
and missing subscripts indicate the matrices of full dimension $s\times s$.
This notation also simplifies the generalized formulation of combined order conditions discussed in \cite{LangSchmitt2022b}.
Multiplying \eqref{OBvStd} by $V_q\T$ from the left and condition \eqref{OBaStd} with a shifted index $n+1\to n$ by $V_r$ from the right yields the necessary conditions
\begin{align}\label{Adqrv}
 \hat A_{q,r}=&\hat B_{q,r}(\sigma)\PP_r^{-1}S_{n,r}^{-1}+\hat K_{q,r}\tilde E_r,\\\label{Adqra}
 \hat A_{q,r}=&\PP_q\T S_{n,q}\hat B_{q,r}(\sigma)-\tilde E_q\T \hat K_{q,r}.
\end{align}
Obviously, $\hat A_{q,r}$ may be eliminated yielding the combined necessary condition
\begin{align}\notag
 \hat K_{q,r}\tilde E_r+\tilde E_q\T \hat K_{q,r}
 =&\PP_q\T S_{n,q}\hat B_{q,r}(\sigma)-\hat B_{q,r}(\sigma)\PP_r^{-1}S_{n,r}^{-1}
 \\\label{LLKd}
 =&\PP_q\T S_{n,q}\QQ_{q,r}(\sigma)\PP_r-\QQ_{q,r}(\sigma)S_{n,r}^{-1}
 =\eins_q\eins\T_r-e_1e_1\T\in\R^{q\times r},
\end{align}
by Lemma~\ref{LBsig} and since $\eins\T\PP_k=\eins\T_k$. We note that the map
\begin{align}\label{LLqr}
\LL_{q,r}:\, X\mapsto X\tilde E_r+\tilde E_q\T X
\end{align}
on the left-hand side of \eqref{LLKd} is singular, and $e_1\T\LL_{q,r}(X)e_1=0$ for any $X$.
For a more detailed discussion see \cite{LangSchmitt2022b,Schmitt2015}.
We like to point out that \eqref{QQform} and \eqref{LLKd} refer to orders $r,q$ used with $\sigma$-uniform order conditions.
The strong restriction \eqref{QQform} for $\QQ$ will limit these orders to $q,r\le s-1$ for sufficiently stable methods.
However, for $r=s-1$ we will try later on to improve the accuracy by applying additional order conditions of order $r_1=s$ for $\sigma=1$ only.
\par
Lemma~\ref{LBsig} not only restricts the form of $\hat B(\sigma)$ to a large extent, it also leaves only a few free parameters in the $\sigma$-independent matrices $A$ and $K$.
\begin{corol}\label{CBAK}
Under the assumptions of Lemma~\ref{LBsig} the coefficient matrices of the standard Peer method satisfy
\begin{align}\label{VBAV}
 \hat B_{q,r}(\sigma)=V_q\T B(\sigma)V_r=e_1\text{\eins}\T\in\R^{q\times r},
 \ \hat A_{q,r}=V_q\T AV_r=\Big(\frac{j-1}{i+j-2}\Big)_{i,j}\in\R^{q\times r},
\end{align}
with $e_1\T\hat Ae_1=\eins\T A\eins=1$.
The matrix $\hat K=V_s\T KV_s$ has Hankel form equaling the Hilbert matrix in its first $q+r-2$ antidiagonals.
\end{corol}
{\bf Proof:}
Multiplying \eqref{QQform} by $\PP_r$ from the right yields the assertion for $V_q\T B(\sigma)V_r$ since $e_1\T\PP_r=\eins\T_r$.
Since the matrix $K$ is diagonal, the matrix $\hat K=\big(\hat\kappa_{ij}\big):=V_s\T KV_s$ and all its submatrices possess Hankel form, $\hat\kappa_{ij}=\xi_{i+j-1},\,1\le i,j\le s$.
Hence, equation \eqref{LLKd} is satisfied as $0=0$ for $i=j=1$ and means
\begin{align}\label{Hankel}
\hat\kappa_{i,j-1}(j-1)+(i-1)\hat\kappa_{i-1,j}=(i+j-2)\xi_{i+j-2}=1
\end{align}
for $3\le i+j\le q+r$.
Non-existent elements of $\hat K$ with index zero are canceled by vanishing factors.
Setting $\ell:=i+j-2\ge 1$, we obtain $\xi_\ell=1/\ell$ for $1\le\ell\le q+r-2$.
\par
Now, equation \eqref{Adqrv} shows that $\hat A_{q,r}=\QQ_{q,r}S_{n,r}^{-1}(\sigma)+\hat K_{q,r}\tilde E_r=e_1e_1\T+\hat K_{q,r}\tilde E_r$, where $\tilde E_r$ shifts columns of $\hat K_{q,r}$ to the right by one.
Hence, the element from the latest anti-diagonal of $\hat K_{q,r}$ entering $\hat A_{q,r}$ is $\hat\kappa_{q,r-1}=\xi_{q+r-2}=1/(q+r-2)$ still belonging to the Hilbert part.
Hence, the formula \eqref{VBAV} follows where the first element of $\hat A_{q,r}$ also fits in by setting $0/0=1$, here.
\qed
\par\noindent
\begin{remark}
The Corollary shows that the diagonal elements of the matrix $K=\diag_i(\kappa_{ii})$ are weights of a quadrature formula associated with the nodes $c_i$, since
\begin{align}\label{Quadform}
 \hat\kappa_{ij}=\sum_{\ell=1}^s\kappa_{\ell\ell}c_\ell^{i+j-2}
  =\frac1{i+j-1}=\int_0^1t^{i+j-2}dt,\ 2\le i+j\le q+r.
\end{align}
Hence, the matrix $K$ is uniquely determined by the nodes for $q+r\ge s+1$.
And if $K$ has only positive diagonal elements it induces a positive quadrature formula where the node polynomial is an orthogonal polynomial for some appropriate weight function, \cite{Xu1994}.
\end{remark}
The design of the standard method will be based on the hat-matrices \eqref{VBAV} since the Corollary leaves few free parameters only in their entries.
However, the matrix $A$ from the standard method has to satisfy a second restriction on its structure since an efficient implementation of the time step requires triangular form of $A$.
This restriction corresponds to a set of $s(s-1)/2$ conditions (6 for $s=4$) which have to be solved with the aid of the remaining parameters in $\hat A$ and the Vandermonde matrix $V$.
With local orders $q=r=s-1$ of the methods the remaining $2s-1$ free parameters (7 for $s=4$) in $\hat A$ may seem to be sufficient for this task without restricting the nodes $c_i$.
Unfortunately, this is not fully true for $s\ge 4$ since linear dependencies lead to some algebraic restrictions on the nodes $c_i$, whose number increases with larger stage numbers $s$.
\begin{lemma}\label{Lnoderest}
The $s$-stage standard method $(A,B(\sigma),K)$ with lower triangular matrix $A$ can satisfy the order conditions \eqref{OBvStd}, \eqref{OBaStd} with $q=r=s-1$ for different stepsize ratios $\sigma_n$ only if the nodes $c_1,\ldots,c_s$ satisfy a set of $(s-2)(s-3)/2$ algebraic equations.
\end{lemma}
{\bf Proof:}
Knowing that $\hat A_{r,r}=V_r\T AV_r$ has fixed entries given by \eqref{VBAV}, for $r=s-1$ we consider the square Vandermonde matrix $\bar V_r\in\R^{r\times r}$ without the last node $c_s$.
Thus, we may write $V_r\T=\bar V_r\T (I_r,v)$ with $v\T=(1,c_s,\ldots,c_s^r)\bar V_r^{-1}$.
Splitting the lower triangular matrix $A$ in the same way, the definition of $\hat A_{r,r}$ is equivalent with
\begin{align*}
\bar V_r\mT \hat A_{r,r}\bar V_r^{-1} =(I_r,v)
 \begin{pmatrix}\bar A&0\\ a_s\T &a_{ss}\end{pmatrix}
 \begin{pmatrix}  I\\v\T \end{pmatrix}
 =\bar A+v(a_s\T+a_{ss}v\T).
\end{align*}
Since the leading block $\bar A$ is also lower triangular, contributions from $A$ above the main diagonal come through the vector $a_s\T+a_{ss}v\T$ only.
Considering now the first $r-1$ entries of the last column number $r$ in this identity, we obtain the linear equation
\begin{align*}
 \big(\bar V_r\mT \hat A_{r,r}\bar V_r^{-1}e_r\big)_{1\ldots r-1}
 =\alpha_r\cdot v_{1\ldots r-1}
\end{align*}
containing one single parameter $\alpha_r:=a_{sr}+a_{ss}v\T e_r$ from $\hat A$ only.
Solutions $\alpha_r$ do only exist if the vectors on both sides are linearly dependent which corresponds to the vanishing of $r-2$ determinants of size $2\times 2$.
This argument repeats for the other columns with shrinking vector length until column number 3 with one $2\times2$-determinant.
This means, one condition for $s=4$, two more for $s=5$, etc.
\qed
\par
The Lemma shows that Peer methods with orders $r=q=s-1$ for variable stepsizes probably do not exist for $s\ge6$ since the number of restrictions on the $s$ nodes is six for $s=6$ and three for $s=5$.

\subsection{Boundary methods}
In addition to the standard method $(A,B(\sigma),K)$ used in the time steps with $1\le n<N$, exceptional methods for the start ($n=0$) and the final step ($n=N$) will be used in order to gain more flexibility in the accurate approximation of the boundary conditions.
And as mentioned before local order higher than two is not possible with triangular coefficients $A_0,A_N,K_0,K_N$.
The order conditions for the starting method \eqref{PeerStrt} are
\begin{align}\label{OBStrtv}
 A_0V_r=&ae_1\T+K_0V_r\tilde E_r,\\
 V_q\T A_0+\tilde E_q\T V_q\T K_0=&\PP_q\T S_{1,q}V_q\T B(\sigma_1)
 =V_q\T A+\tilde E_q\T V_q\T K,
\end{align}
where the forward condition is copied from \cite{LangSchmitt2022a} and the adjoint condition uses \eqref{OBaStd} twice with $n=0$ and $n>0$.
Hence, both conditions are independent of the ratio $\sigma_1$ and also $A_0,K_0$ may be constant.
\par
As a change to earlier papers \cite{LangSchmitt2022a,LangSchmitt2022b}, we will now use the  matrix $B(\sigma)$ from the standard method also in the end step $n=N$.
Hence, for the end method, the forward condition \eqref{OBvsw} is complemented with only one adjoint condition yielding
\begin{align}
  A_NV_r-K_NV_r\tilde E_r=&B(1)V_r\PP_r^{-1},\\\label{OBaRB}
  A_N\T V_q+K_N\T V_q\tilde E_q=&w\eins\T,
\end{align}
where the first condition uses \eqref{OBvsw} and second condition comes from \cite{LangSchmitt2022a}.
Since $\tilde E_re_1=0$, first consequences of \eqref{OBStrtv}, \eqref{OBaRB} are
\begin{align}\label{aundw}
  a=A_0\eins,\quad w=A_N\T\eins.
\end{align}
Order $r$ of accuracy for the extrapolation step \eqref{PeerEnd} also requires that, see \cite{LangSchmitt2022b},
\begin{align}\label{EExtrp}
 w\T V_r=\eins_r\T.
\end{align}
Comparing with \eqref{aundw}, this means that $\eins\T A_NV_r=e_1\T V_s\T AV_r=\eins_r\T$ and hence the $r$ first entires of $\hat A_N$ coincide with those from $\hat A$, see \eqref{VBAV}.
\par
Similar to \eqref{LLKd}, the order conditions for each boundary method can be combined to the following necessary conditions.
They are much simpler now than those in \cite{LangSchmitt2022b} due to Lemma~\ref{LBsig} which showed that $\QQ_{q,r}=e_1e_1\T$.
These combined conditions are
\begin{align}\label{KOBstrt}
 \LL_{q,r}(V_q\T K_0 V_r)=\PP_q\T\QQ_{q,r}\PP_r-V_q\T ae_1\T
 =&\eins\eins\T-V_q\T a e_1\T,
\\\label{KOBend}
 \LL_{q,r}(V_q\T K_NV_r)=\eins\eins\T-\QQ_{q,r}
 =&\eins\eins\T-e_1e_1\T,
\end{align}
with the map $\LL_{q,r}$ defined in \eqref{LLqr}.
Here, the singularity of this map with the property $e_1\T \LL_{q,r}e_1=0$ requires that $e_1\T \QQ_{q,r}e_1=1$ in \eqref{KOBend} which was the reason for assuming the normalization \eqref{eAe} above.
And in \eqref{KOBstrt} it requires that $1=e_1\T V_q\T a=\eins\T a=\eins\T A_0\eins$.
\par
An additional constraint comes from the one-leg-condition.
Rewriting \eqref{OBU} with $q_b=2$ for $\hat K_0$, we see that $\cc\T K_0V_r=e_2\T V_q\T K_0V_r=e_1\T K_0CV_r$, where the product $CV_r$ causes a shift in the columns of $V_{r+1}$ to the left.
This means that $\hat K_0$ has Hankel form in its  first two rows, $\hat\kappa_{1,j}^{(0)}=\hat\kappa_{2,j-1}^{(0)},\,j\ge 2$.
And similar to \eqref{Hankel}, this property carries over to $\LL_{2,r}(V_2\T K_0V_r)$.
Hence, \eqref{KOBstrt} has solutions satisfying \eqref{OBU} only if $V_{2}\T a=e_1\in\R^{2}$, which by \eqref{aundw} means
\begin{align}\label{A0lsb}
 \eins\T A_0\eins=1,\quad \cc\T A_0\eins=0.
\end{align}
Summing up, the three methods in the Peer triplet use a set of seven coefficient matrices $(A_0,K_0)$, $(A,B(\sigma),K)$, $(A_N,B(\sigma_N),K_N)$, where the parameter-dependent matrix $B(\sigma)$ is the same in all time steps.
The combined equations \eqref{KOBstrt}, \eqref{KOBend} allow for a detached design of the standard method $(A,B(\sigma),K)$ without reference to the boundary methods since all necessary restrictions for their existence are known.

\section{Boundedness of matrix families}\label{SBMF}
A severe new difficulty  arising for variable-stepsize Peer methods is uniform zero-stability which requires that arbitrary long products
\begin{align}\label{Bprods}
 \bar B_n\bar B_{n-1}\cdots \bar B_k,\quad 0\le k\le n,
\end{align}
of stability matrices $\bar B_n:=A^{-1}B_n=A^{-1}B(\sigma_n)$ are uniformly bounded.
We note that the stability matrices of the adjoint steps \eqref{AdjABK} are related by $A\mT B(\sigma)\T=\big(B(\sigma)A^{-1})\T$.
Uniform boundedness of long products is connected with the problem of the {\em joint spectral radius} of matrix families.
However, to our knowledge no general, applicable theory exists for this problem, see, e.g. \cite{BlondelNesterov2010}.
Instead, as in \cite{SchmittWeiner2017}, we incorporate the construction of a suitable norm satisfying
$\ltnorm\bar B_n\rtnorm=1$ for all interior steps $n$ into the design and search process of the standard peer method.
\par
As a consequence of the conditions \eqref{precons}, the matrix  $\bar B_n$ possesses the eigenvalue one with known eigenvectors and the super-convergence effect that will be exploited later on requires that one is an isolated dominant eigenvalue whose absolute value exceeds that of all others sufficiently, see \eqref{lbd2}.
Hence, our objective is to construct a constant similarity transformation which separates the eigenvalue one in block diagonal form such that a simple norm of the lower block is uniformly bounded below one for stepsize ratios $\sigma_n$ from a reasonably large interval $[\underline\sigma,\bar\sigma]\ni 1$.
\par
A first clue to the construction of such a norm gives Corollary~\ref{CBAK} which we discuss in detail here for the case of interest with $q=r=3<s=4$. The Corollary shows that
\begin{align}\label{BdAd4}
  \hat A=\begin{pmatrix} 1&1&1&\hat a_{14}\\
  0&\frac12&\frac23&\hat a_{24}\\
  0&\frac13&\frac12&\hat a_{34}\\
 \hat a_{41}&\hat a_{42}&\hat a_{43}&\hat a_{44}
 \end{pmatrix},\quad
\hat B(\sigma)=\begin{pmatrix} 1&1&1&\hat a_{14}\\
  0&0&0&\hat b_{24}(\sigma)\\
  0&0&0&\hat b_{34}(\sigma)\\
 \hat a_{41}&\hat b_{42}(\sigma)&\hat b_{43}(\sigma)&\hat b_{44}(\sigma)
 \end{pmatrix},
\end{align}
and $\hat K=hankel(1,\frac12,\frac13,\frac14,\hat k_5,\hat k_6,\hat k_7)$.
We confirm that the first rows and columns of $\hat A,\hat B$ are identical and independent of $\sigma$ due to \eqref{precons}.
The elements in their fourth rows and columns may be subject to further restrictions.
Written with these transformed matrices $\hat A,\hat B,\hat K$, the full set of order conditions \eqref{OBvStd},\eqref{OBaStd} simplifies.
With selection matrices $\hat I_k={I_k\choose 0\T}\in\R^{s\times k}$ it reads
\begin{align}\label{OBDv}
 \Big(\hat A-\hat K\tilde E_s-\hat B(\sigma) P_s^{-1}S^{-1}\Big)\hat I_r=&\,0,\\\label{OBDa}
 \hat I_q\T\Big(P_s\mT(\hat A+\tilde E_s\T \hat K)-S\hat B(\sigma)\Big)=&\,0,
\end{align}
and may be solved easily for $q,r\ge3$ by
\begin{align}\label{Bsigma}
 \begin{array}{rl}
 \hat b_{42}=&\hat a_{41}+\sigma\hat a_{42}-\frac{1}4\sigma,\\
 \hat b_{43}=&\hat a_{41}+\sigma(2\hat a_{42}-\frac12)+\sigma^2(\hat a_{43}-2\hat k_5),\\
 \hat b_{24}=&\sigma^{-1}\big(\frac14-\hat a_{14}+\hat a_{24}\big),\\
 \hat b_{34}=&\sigma^{-2}\big(-\frac12+\hat a_{14}-2\hat a_{24}+\hat a_{34}+2\hat k_5\big).
 \end{array}
\end{align}
For $q=r=3$, only the last element $\hat b_{44}(\sigma)$ of $\hat B(\sigma)$ is an unrestricted parameter.
And we remind that other elements except $\hat a_{14},\hat a_{41}$ are determined by requiring triangularity of the matrix $A=V\mT\hat AV^{-1}$.
\par
The stability matrix $\bar B(\sigma)=A^{-1}B(\sigma)$ of the standard method is similar to the matrix
\begin{align}\label{hatM}
 \hat A^{-1}\hat B(\sigma)=V_s^{-1}A^{-1}B(\sigma)V_s=V_s^{-1}\bar B(\sigma)V_s,
\end{align}
and an analogous result holds for the adjoint stability matrix $(B(\sigma)A^{-1})\T$.
Since the factor $L_A$ in the LU-decomposition $\hat A=L_AU_A$ inherits the first column from $\hat A$ and $U_A$ its first row, a further stepsize-independent similarity yields a very simple structure for both stability matrices through
\begin{align}\label{LBR}
 B_{LU}(\sigma)
 :=&L_A^{-1}\hat B(\sigma)U_A^{-1}
 =U_A\big(\hat A^{-1}\hat B(\sigma)\big)U_A^{-1}
 =L_A^{-1}\big(\hat B(\sigma)\hat A^{-1}\big)L_A
 \\\notag
 =&\left(\begin{array}{c|c}
 1&\\\hline
 &B_{se}(\sigma)
 \end{array}\right)
 =\left(\begin{array}{c|ccc}
  1&0&0&0\\\hline
  0&0&0&*\\
  0&0&0&*\\
  0&*&*&*
 \end{array}\right)
\end{align}
where the last matrix is a sketch for $s=4$ and $q=r\ge3$.
Obviously, $B_{LU}(\sigma)$ has block diagonal form separating the dominant eigenvalue one from the rest,
where the asterisks represent elements of $B_{se}(\sigma)$ which may be nontrivial.
The tailored norm $\ltnorm\bar B_n\rtnorm=1$ will use the maximum norm for the southeast block $B_{se}(\sigma)$ with an additional scaling.
\par
As in \cite{SchmittWeiner2017}, we look for additional scalings with a diagonal matrix 
$\bar\Omega=\diag(\omega_2,\ldots,\omega_s)$ such that
\begin{align}\label{UBlkNrm}
 \|\bar\Omega^{-1} B_{se}(\sigma)\bar\Omega\|_\infty\le \tilde\gamma<1,
\end{align}
for an appropriate set of stepsize ratios $\sigma$.
Then, for the positive vector $\bar\omega:=\bar\Omega\eins>0$, the inequality 
$|B_{se}(\sigma)|\bar\omega\le\tilde\gamma\bar\omega$ holds and quasi-optimal scalings may be computed efficiently by solving a linear program:
\begin{align}\label{LinProgSkal}
 \min\eins\T \bar\omega:\ \bar\omega\ge\eins,\ \big(\tilde\gamma I-|B_{se}(\sigma)|\big)\bar\omega\ge0\ \text{ for all }\sigma\in\Sc,
\end{align}
with a small set of stepsize ratios, e.g., $\Sc:=\{1/\hat\sigma,1,\hat\sigma\},\,\hat\sigma>1$.
Since $B_{se}(\sigma)$ is a continuous function, we found that the norm bound $\|\bar\Omega^{-1}B_{se}(\sigma)\bar\Omega\|_\infty\le1$ may hold in rather large $\sigma$-intervals $[\underline\sigma,\bar\sigma]\supseteq\{1/\hat\sigma,1,\hat\sigma\}$ by appropriate definitions of the free parameter $\hat b_{44}(\sigma)$.
\begin{remark}
By \eqref{Bsigma} the element functions in the last row and column of $\hat B(\sigma)$ are determined by order conditions.
The only exception is its last element which may be chosen freely for $q,r\le s-1$ as an arbitrary function $\hat b_{ss}(\sigma)$.
Now, due to the triangular forms of $L_A$ and $U_A$ the element $\hat b_{ss}(\sigma)$ appears in the last element of $B_{se}$ only having the form $b_{ss}^{LU}(\sigma)=b_{ss}^\ast(\sigma)+\hat b_{ss}(\sigma)/u_{ss}$, $u_{ss}=e_s\T U_Ae_s$.
The part $b_{ss}^\ast(\sigma)$ which does not depend on $\hat b_{ss}(\sigma)$ is easily computed from $L_A,U_A$ and the rest of $\hat B$.
The diagonal element $b_{ss}^{LU}(\sigma)$ is not affected by the weight $\bar\Omega$ in the norm $\|\bar\Omega^{-1}B_{se}(\sigma)\bar\Omega\|_\infty$ and contributes to it in the last entry only as $|b_{ss}^{LU}(\sigma)|= | b_{ss}^\ast(\sigma)+\hat b_{ss}(\sigma)/u_{ss}|$.
Now, this term may be canceled exactly by choosing
\begin{align}\label{minbt44}
 \hat b_{ss}(\sigma)=-u_{ss} b_{ss}^\ast(\sigma).
\end{align}
The computation of the diagonal scaling $\bar\Omega$ may then follow up by solving \eqref{LinProgSkal}.
Hence, the transformation \eqref{LBR} allows for an explicit  norm-optimal choice of the free parameter function $\hat b_{ss}(\sigma)$ in $\hat B(\sigma)$.
This is helpful in some cases where the matrix norm is the bottleneck.
However, the choice of $\hat b_{ss}(1)$ also influences $A(\alpha)$-stability and other criteria which may be more important in other cases.
\end{remark}
Summarizing this discussion, we will verify the assumption of uniform zero-stability in the convergence theorem below for each standard method $(A,B(\sigma),K)$ by presenting explicitly some weight matrix $W\in\R^{s\times s}$ and an interval for which holds
\begin{align}\label{GlmStab}
  \ltnorm\bar B(\sigma)\rtnorm:=
  \|W^{-1} \bar B(\sigma)W\|_\infty=1,\ \sigma\in[\underline\sigma,\bar\sigma],\, 0<\underline\sigma<1<\bar\sigma.
\end{align}
Since all involved matrix elements are known with rational numbers or square roots, this property may even be proven rigorously.
According to \eqref{LBR}, the weight matrix may have the form $W=V_s U_A^{-1}\Omega$ with the extended weight $\Omega=\diag(1,\omega_2,\ldots,\omega_s)$.
A simple property of this matrix is $We_1=\eins$.
For the methods designed later on, we may approximate other entries of $W$ with simpler rational expressions and present these more usable weight matrices explicitly.
Since $W^{-1}\bar B(\sigma)W=W^{-1}A^{-1}B(\sigma)W=(AW)^{-1}B(\sigma)A^{-1}(AW)$, a similar norm for the adjoint stability matrix $\tilde B(\sigma_n)\T= \big(B(\sigma_n)A^{-1}\big)\T$ is given with the matrix $W^\dagger=(AW)\mT$ satisfying
\begin{align}\label{Stbsig}
 \ltnorm\tilde B(\sigma)\T\rtnorm:=\|(W^\dagger)^{-1} A\mT B(\sigma)\T W^\dagger\|_1=1,\ \sigma\in[\underline\sigma,\bar\sigma].
\end{align}
The super-convergence effect which prevents the usual loss of one order between local and global error requires that long products of stability matrices converge to the following rank-1-matrices
\begin{align}\label{Binfty}
 \bar B^\infty=\eins\eins\T A,\quad (\tilde B^\infty)\T=\eins\eins\T A\T,
\end{align}
i.e. $ \bar B_n\bar B_{n-1}\cdots \bar B_k\to\bar B^\infty\ (n\to\infty)$, and $\tilde B_n\T\cdots\tilde B_k\T\to(\tilde B^\infty)\T$ $(k\to\infty)$.
These limits are composed of the left and right eigenvectors to the dominant eigenvalue which are both known due to \eqref{precons}.
And the dominant eigenvalue in \eqref{Binfty} is $\eins\T A\eins=e_1\T \hat Ae_1=1$ indeed, see \eqref{eAe}.
It is essential that the weight matrices $W,W^\dagger$ based on the transformation to block diagonal form \eqref{LBR} are independent of $\sigma_n$.
Hence, this structure carries over to arbitrarily long products
\begin{align}\label{BProd}
 B_{LU}(\sigma_n)\cdots B_{LU}(\sigma_k)=
 \begin{pmatrix}
  1&0\\0& B_{se}(\sigma_n)\cdots B_{se}(\sigma_k)
 \end{pmatrix}.
\end{align}
In fact, super-convergence requires that the lower block in \eqref{BProd} tends to zero for $n-k\to\infty$ for admissible stepsize ratios $\sigma_n,\ldots,\sigma_k$ and that infinite sums of products of $B_{se}$-matrices are uniformly bounded.
If \eqref{UBlkNrm} is true for $\sigma\in[\underline\sigma,\bar\sigma]$, then the following norm estimates hold:
\begin{align}\label{NormSupKv}
 \ltnorm\bar B(\sigma)-\bar B^\infty\rtnorm
  =&\;\|W^{-1}\big(\bar B(\sigma)-\bar B^\infty\big)W\|_\infty
  =\|W^{-1}\bar B(\sigma)W-e_1e_1\T\|_\infty\le\tilde\gamma<1,\\\label{NormSupKa}
 \ltnorm\tilde B(\sigma)\T-(\tilde B^\infty)\T\rtnorm
  =&\;\|(W^\dagger)^{-1}\big(\tilde B(\sigma)\T-(\tilde B^\infty)\T\big)W^\dagger\|_1
  =\|(W^\dagger)^{-1}\tilde B(\sigma)\T W^\dagger-e_1e_1\T\|_1\le\tilde\gamma,
\end{align}
since $We_1=\eins$ and $\eins\T AW=e_1\T\hat AU_A^{-1}\bar \Omega=e_1\T$.
We note that since $\bar B^\infty$ and $(\tilde B^\infty)\T$ share both their dominant invariant subspaces with all matrices $\bar B_n$, resp. $\tilde B_n\T$ (maybe except the boundary steps) these bounds transfer to products, e.g., for $1\le n< k<N$,
\begin{align}\notag
  \ltnorm\tilde B(\sigma_n)\T\cdots\tilde B(\sigma_{k-1})\T-(\tilde B^\infty)\T\rtnorm
  =&\;\ltnorm (\tilde B(\sigma_n)\T-(\tilde B^\infty)\T)\cdots(\tilde B(\sigma_{k-1})\T-(\tilde B^\infty)\T)\rtnorm\\
  \label{Btprest}
  \le&\;\tilde\gamma^{k-n}.
\end{align}

\section{Design of 4-stage Peer triplets for variable stepsizes}\label{SDesign4}
Efficient Peer triplets need to satisfy more requirements than the order conditions and stability properties discussed so far.
Before starting with the actual construction of methods, we mention a few additional useful properties.
Time integration methods usually loose one order between local and global error.
This loss may be prevented through some super-convergence effect for the standard method, see \cite{LangSchmitt2022b}.
Surprisingly, in our present setting the corresponding additional order conditions are fulfilled automatically for $\sigma=1$.
The lemma refers to the error vectors $\beta_r,\beta_q^\dagger$ from \eqref{lokFRestv}, \eqref{lokFResta}.
\begin{lemma}\label{LSupKBed}
Let the standard method $\big(A,B(\sigma),K\big)$ satisfy the order conditions \eqref{OBvStd}, \eqref{OBaStd} with $q=r=s-1$ for two different stepsize ratios $\sigma_n$, at least.
\par\noindent
a) Then, also
\begin{align}\label{SupKvv}
 r!\eins\T A\beta_{r}(\sigma)=\eins\T\big(A\cc^r-B(\sigma)(\cc-\eins)^r\sigma^{-r}-r K\cc^{r-1}\big)&=(\hat a_{1s}-1)(1-\sigma^{-r}),\\\label{SupKva}
 q!\eins\T A\T\beta_{q}^\dagger(\sigma)=\eins\T\big(A\T\cc^q-B(\sigma)\T(\eins+\sigma\cc)^q+qK\cc^{q-1}\big)
 &=\hat a_{s1}(1-\sigma^q),
\end{align}
where $(\hat a_{ij})=\hat A=V_s\T AV_s$.
\par\noindent
b) Let the standard method also satisfy $\hat a_{1s}=1$ (i.e. $e_1\T\hat A=\eins\T$) and $\hat k_{s}=1/s$, then
\begin{align}\label{SupKvs}
 s!\eins\T A\beta_s(\sigma)=\eins\T\big(A\cc^s-B(\sigma)(\cc-\eins)^s\sigma^{-s}-s K\cc^{s-1}\big)&=(1-\sigma^{-s})(\eins\T A\cc^s-1).
\end{align}
\end{lemma}
{\bf Proof:}
a) From Corollary~\ref{CBAK} and condition \eqref{precons} we know that
$e_1\T\hat B(\sigma)=(1,\ldots,1,\hat a_{1s})=\eins\T+(\hat a_{1s}-1)e_s\T$ is independent of $\sigma$, and $e_1\T\hat Ke_{s-1}=1/(s-1)$.
In analogy to \eqref{OBDv}, we see that the first assertion is equivalent with
\begin{align*}
 &\;e_1\T\big(\hat Ae_s-\hat B(\sigma)\PP_s^{-1}e_s\sigma^{-r}-(s-1)\hat Ke_{s-1}\big)\\
 =&\; \hat a_{1s}-(\eins\T+(\hat a_{1s}-1)e_s)\PP_s^{-1}e_s\sigma^{-r}-1
 =(\hat a_{1s}-1)(1-\sigma^{-r}),
\end{align*}
since $\eins\T\PP_s^{-1}=e_1\T$.
Now, the first column of $\hat B$ is $\hat B(1)e_1=e_1+\hat a_{s1}e_s$ and rewriting \eqref{SupKva} as before we obtain
\begin{align*}
 &(e_s\T\hat A-e_s\T\PP_s\T S_s\hat B(\sigma)+(s-1)e_{s-1}\T\hat K\big)e_1
 = \hat a_{s1}-e_s\T\PP_s\T(e_1+\sigma^q\hat a_{s1}e_s)+1=\hat a_{s1}(1-\sigma^q).
\end{align*}
b) Due to $\eins\T B(\sigma)\equiv\eins\T A$, we have for $r=s$ that
\begin{align*}
& \eins\T \big(A\cc^s-B(\sigma)(\cc-\eins)^s\sigma^{-s}-s K\cc^{s-1}\big)\\
&=(1-\sigma^{-s})\eins\T A\cc^s-\sum_{k=0}^{s-1}{s\choose k}(-1)^k\eins\T B(\sigma)\cc^k\sigma^{-s}-s\hat\kappa_{1,s}=(1-\sigma^{-s})(\eins\T A\cc^s-1).
\end{align*}
The last step used $\eins\T B\cc^k=e_1\T\hat B e_{k+1}=e_1\T\hat Ae_{k+1}=1$, $k\le s-1$.
\qed
\par
We note that the expression \eqref{SupKvv} vanishes under the assumptions of part b) and that the assumption $\hat\kappa_{1s}=1/s$ is automatically satisfied for $s\ge 4$ due to Corollary~\ref{CBAK}, where $q+r-2=2s-4\ge s$.
\par
The conditions \eqref{SupKvv}, \eqref{SupKva} will cancel the leading terms $\bar B^\infty\beta_r(1),\,(\tilde B^\infty)\T \beta_q^\dagger(1)$ in the global error.
However, it is also required that products of stability matrices like \eqref{Bprods} converge to the rank-1-matrix $\bar B^\infty$ in \eqref{Binfty}.
The norm bounds \eqref{NormSupKv}, \eqref{NormSupKa} ensure this in theory and they will be used to verify an assumption in the main convergence theorem with an upper bound $\tilde\gamma$ quite close to one.
However, in computational practice, convergence to $\bar B^\infty$ has to be fast enough.
So, according to \cite{LangSchmitt2022b} a similar bound with a smaller limit will also be enforced, but only for the eigenvalues.
The requirement
\begin{align}\label{lbd2}
 |\lambda_2\big(A^{-1}B(1)\big)|\le \gamma<1\quad\text{ with } \gamma\approx 0.8
\end{align}
for the absolutely second largest eigenvalue works well in practice.
The conditions from Lemma~\ref{LSupKBed} are inner products with the error vectors canceling only the leading error terms.
In order to cover their full impact, we also monitor the whole error vectors in the design process and introduce
\begin{align}\label{FeKov}
err_{r,n}=&\frac{1}{r!}\|\cc^r-A_n^{-1}B(1)(\cc-\eins)^r-rA_n^{-1}K_n\cc^{r-1}\|_{\infty}, \\\label{FeKoa}
err_{q,n}^\dagger =&\frac{1}{q!}\|\cc^q -A_n\mT B(1)\T(\eins+\cc)^q+qA_n\mT K_n\cc^{q-1}\|_{\infty},
\end{align}
again as the essential error constants \cite{LangSchmitt2022b}.
Modifications are required at the boundaries as
\begin{align*}
 err_{r,0}=\frac1{r!}\|\cc^r-rA_0^{-1}K_0\cc^{r-1}\|_{\infty},\quad
 err_{q,N}^\dagger=\frac1{q!}\|\cc^q+A_N\mT(q K_N\T\cc^{q-1}-w)\|_{\infty}.
\end{align*}
For the standard method, $1\le n<N$, the additional index $n$ is omitted, as usual.
Solving the nonlinear equations \eqref{PeerStd}, \eqref{PeerStrt} requires non-singularity of the re-arranged Jacobians
$K_n^{-1}A_n-h_n\nabla_Y F$.
In the context of stiff problems an appropriate assumption is that
\begin{align}\label{defmun}
 \mu_n:=\min_j\mbox{Re}\,\lambda_j(K_n^{-1}A_n)>0,\;n=0,\ldots,N,
\end{align}
for the eigenvalues of $K_n^{-1}A_n$.
In the design process of the boundary methods further data are monitored like the size of the norms $\|A_N^{-1}B(1)\|,\ \|A_0\mT B(1)\T\|$ or $K_0,K_N$.
\par
Considering the assertions of Lemma~\ref{LSupKBed}, we see that the choices of the method parameters lead to different consequences in \eqref{SupKvv}, \eqref{SupKva}.
On the one hand, for $\hat a_{1s}=1$ resp. $\hat a_{s1}=0$ the conditions for super-convergence are satisfied for arbitrary stepsize ratios $\sigma_n$, i.e. for arbitrary grids.
However, for  $\hat a_{1s}\not=1$ resp. $\hat a_{s1}\not=0$ both expressions have the form $O(1-\sigma_n)=O(h_n)$ in the $n$-th time step if we consider smooth grids satisfying $\sigma_n=1+O(h_n)$ only.
Hence, different design decisions lead to methods with different properties and scopes of application.
In the search process it was observed that very favorable properties of the standard method coincide with the choice $c_4=1$ for the forward time step, and to a lesser extent with $c_1=0$ for the adjoint step.
Higher-order methods of this form have a surprising property which is not obvious in our redundant formulation \eqref{PeerStd} of the Peer method.
\begin{lemma}\label{LLSRK}
Let the standard method $(A,B,K)$ satisfy the order conditions \eqref{OBvsw}, \eqref{OBasw} $\sigma$-uniformly for $r=q=s-1$.
\par\noindent a) If also $e_1\T\hat A=\eins\T$ and $c_s=1$ hold, then
\begin{align}\label{LSRKBv}
 \eins\T A=e_s\T,\quad e_s\T\bar B(\sigma)=e_s\T A^{-1}B(\sigma)=e_s\T.
\end{align}
b) If also $\hat Ae_1=e_1$ and $c_1=0$ is satisfied, then
\begin{align}\label{LSRKBa}
 A\eins=e_1,\quad e_1\T\tilde B(\sigma)=e_1\T A\mT B(\sigma)\T=e_1\T
\end{align}
\end{lemma}
{\bf Proof:}
a) Assumption $c_s=1$ means that $e_s\T V=\eins\T$.
Hence,
\[ \eins\T A=e_1\T \hat A V^{-1}=\eins\T V^{-1}=e_s\T.\]
Now, \eqref{precons} also gives $\eins\T B(\sigma)=\eins\T A=e_s\T$ and shows that the last row of the stability matrix satisfies
$e_s\T A^{-1}B(\sigma)=\eins\T B(\sigma)=e_s\T$.
\par\noindent
b) The assumption means $V\mT\hat Ae_1=A\eins=V\mT e_1=e_1$ since $e_1\T V=e_1\T$ for $c_1=0$.
Hence, by \eqref{precons} follows $B(\sigma)A^{-1}e_1=B(\sigma)\eins=A\eins=e_1$.
 \qed
\begin{remark}
Part a) of the Lemma shows that $A$ has vanishing column sums with exception of the last one and that $e_s\T$ is the dominant left eigenvector of the stability matrix $\bar B(\sigma)$.
Multiplying the time step \eqref{PeerStd} by $A^{-1}$ we see that now the final stage of the standard Peer method reads
\begin{align}\label{LSRK}
  Y_{ns}=Y_{n-1,s}+h_n\sum_{j=1}^s\kappa_{jj}g(Y_{nj},P_{nj}),
\end{align}
where $Y_{n-1,s}\cong y(t_n)$ since $c_s=1$.
Obviously, this last stage corresponds to the final one of a Runge-Kutta method and we call this form LSRK {\em(Last Stage is Runge-Kutta)}.
In fact, the quadrature weights of this stage do not depend on $\sigma$ since $e_s\T A^{-1}K=\eins\T K$, see \eqref{Quadform}.
Under the assumptions of part b) of the Lemma stage number one of the adjoint time step, being the last step backwards, also is a final Runge-Kutta stage of the form
\begin{align}\label{LSRKa}
  P_{n1}=P_{n+1,1}-h_n\sum_{j=1}^s\kappa_{jj}\varphi(Y_{nj},P_{nj}),
\end{align}
since again $e_1\T A\mT K=\eins\T K$.
\par
Compared to other candidates, the properties of such methods are surprisingly favorable since the LSRK form seems to decouple the two time steps to some extent which may weaken the difficulties brought by the $\sigma$-dependence of $\bar B(\sigma)$.
Although \eqref{LSRK}, \eqref{LSRKa} look like simple Runge-Kutta end steps we  stress the fact that the Peer method avoids order reduction by using the increments $g(Y_{nj},U_{nj})$, $\varphi(Y_{nj},U_{nj})),\,1\le j\le s,$ having high stage order due to their two-step structure.
\end{remark}
It is a general observation in the design of numerical integrators that there may be a trade-off between accuracy requirements and stability properties.
This happens also for Peer triplets and we will have to discuss different design choices.
We remind of the general remark at the end of Section~\ref{SAKconst} that the knowledge of all necessary restriction for the boundary methods allows a detached development of the standard method having a small number of free parameters only.
The choice of matching boundary methods may then follow up.

\subsection{The triplet \texttt{AP4o33vg} for general grids}
First, we discuss a method with orders $r=q=s-1=3$ where both conditions \eqref{SupKvv} and \eqref{SupKva} for super-convergence are cancelled $\sigma$-uniformly.
Hence, there will be no smoothness restriction on the grids, only the usual bounds on the stepsize ratios, $\sigma_n\in[\underline{\sigma},\bar\sigma]$, required by uniform zero stability.
In our naming scheme from previous papers \cite{LangSchmitt2022a,LangSchmitt2022b}, the triplet is denoted as \texttt{AP4o33vg} ({\bf A}djoint {\bf P}eer method with {\bf 4} stages and {\bf o}rders {\bf 3,3} for {\bf v}ariable stepsizes on {\bf g}eneral grids).
By Lemma~\ref{LBsig} the necessary conditions for solvability of the combined conditions \eqref{KOBstrt}, \eqref{KOBend} are fulfilled and a separate construction of the standard method is possible.
Considering such Peer methods, the largest angles of $A(\alpha)$-stability were observed with $c_1=0$ and $c_s=1$.
By Lemma~\ref{LLSRK} this means that the method has LSRK form forward and backward in time.
In this situation, only 2 free parameters remain in \eqref{BdAd4}, the node $c_2$ and $\hat b_{44}(\sigma)$ since the remaining entries of $\hat A$ as well as $c_3$ are restricted by the required triangularity of $A$ itself.
The largest angle $\alpha$ was seen for $c_2=\frac13$ yielding a method of very simple form with
\begin{align}\label{K4o33vg}
 \cc\T=\left(0,\frac13,\frac23,1\right),\quad
  K=\diag\left(\frac18,\frac38,\frac38,\frac18\right),\quad
 A=\begin{pmatrix}
  1&0&0&0\\[1mm]
  -\frac94&\frac94&0&0\\[1mm]
  \frac94&-\frac92&\frac94&0\\[1mm]
  -1&\frac94&-\frac94&1
 \end{pmatrix}.
\end{align}
We note that the positive definite matrix $K$ contains the weights of the {\em pulcherrima} quadrature rule associated with the given nodes in $\cc$.
With a fine-tuned choice of $\hat b_{44}(\sigma)$ the method possesses a stability angle of $\alpha\doteq 61.59$ and an interval of uniform zero stability $[\underline{\sigma},\bar\sigma]=[0.57,1.75]$ associated with the weight matrix $W$ given in Appendix~A1.
The damping factor in \eqref{lbd2} is $\gamma=0.8$ and the error constants coincide $err_3=err_3^\dagger\doteq 0.0093$.
For later reference, we display the essential properties of the local errors in \eqref{lokFRestv}, \eqref{lokFResta} with $r=q=3=s-1$ here:
\begin{align}\label{LF4o33vg}
 \begin{array}{rlr}
 \tau_n^Y=&h_n^r\beta_r(\sigma_n)y^{(r)}(t_n)+O(h_n^s\|y^{(s)}(t_n)\|_{[n]}),
  &\eins\T A\beta_r(\sigma)\equiv0\text{ for all }\sigma,\\[2mm]
 \tau_n^P=&h_n^q\beta_q^\dagger(\sigma_{n+1})p^{(q)}(t_n)+O(h_n^s\|p^{(s)}\|_{[n]}),
 & \eins\T A\T\beta_q^\dagger(\sigma)\equiv0\text{ for all }\sigma.
 \end{array}
\end{align}
A remarkable property of this method is that it coincides with its own adjoint, similar to the situation for BDF methods in \cite{LangSchmitt2022a,SchroederLangWeiner2014}.
In fact, with the flip permutation $\Pi=\big(\delta_{i,s+1-i}\big)$ it holds that $\Pi A\Pi=A\T$ and $\Pi B(1)\Pi=B(1)\T$ leading to the identity $\Pi A^{-1}B(1)\Pi=A\mT B(1)\T$ for the two stability matrices.
As mentioned before, order 3 in the boundary steps is not possible with triangular matrices $K_n,A_n$, $n=0,N$.
In order to extend the symmetry of the standard method as far as possible, we kept the weight matrices $K_n\equiv K,\,n=0,\ldots,N$, and considered rank-1-perturbations for $A_0$ and $A_N$ in block triangular form satisfying the order conditions for these methods.
A very attractive choice is
\begin{align}\label{RM4o33vg}
 A_0=
 \begin{pmatrix}
  \frac{49}{80} & \frac34 & -\frac3{16}& 0\\[1mm]
  -\frac{87}{80}& 0 & \frac9{16}& 0\\[1mm]
  \frac{87}{80} & -\frac94 & \frac{27}{16}& 0\\[1mm]
  -\frac{49}{80}& \frac32 & -\frac{33}{16}& 1
 \end{pmatrix},\quad
 A_N=\begin{pmatrix}
  1&  0&  0&  0\\[1mm]
  -\frac{33}{16} & \frac{27}{16}& \frac9{16}& -\frac3{16}\\[1mm]
  \frac32 & -\frac94 & 0& \frac34\\[1mm]
  -\frac{49}{80}& \frac{87}{80}& -\frac{87}{80}&  \frac{49}{80}
 \end{pmatrix}
 =(\Pi A_0\Pi)\T,
\end{align}
satisfying \eqref{defmun} with $\mu_0=2.7$.
Since even the vectors from \eqref{PeerStrt} and \eqref{PeerEnd} are related by $w=\Pi a$ we see that the whole triplet coincides with its own adjoint.
This is a very nice property and we will occasionally address this  triplet \texttt{AP4p33vg} as {\em pulcherrima}.
However, this beauty comes with a price, there are other triplets with much larger stability angles.
\par
An overview of all standard methods and their essential properties is presented in Table~\ref{TPT} which displays the names of the triplets, the number of stages $s$, the orders $(r,q)$ satisfied for all $\sigma$ and the order $r_1$ satisfied for $\sigma=1$ only, the range of the nodes, the stability angle $\alpha$, the interval of stepsize ratios $[\underline\sigma,\bar\sigma]$ for which strong zero-stability \eqref{Stbsig} holds and finally the forward and adjoint error constants.
Table~\ref{TPB} contains properties of the boundary methods, for instance the block form \texttt{blksz} of the matrices $K_n^{-1}A_n$, $n=0,N$, which defines the size of the stage systems to be solved in these 2 time steps.

\subsection{The triplet \texttt{AP4o33vs} for smooth grids}
The LSRK property \eqref{LSRKBv} for the state equation alone already leads to improved properties of the standard method.
We now consider such methods with $r=q=s-1=3$ where only the condition \eqref{SupKvv} is canceled $\sigma$-uniformly while the adjoint condition \eqref{SupKva} is satisfied for $\sigma=1$ only with the benefit of having 2 additional free parameters, $\hat a_{s1}$ and $c_1$.
\par
Monte-Carlo-type computer searches, incorporating the construction of a suitable norm $\ltnorm\cdot\rtnorm$ in \eqref{GlmStab} in Section~\ref{SBMF}, found standard methods with stability angles of over 84 degrees.
Restrictions on the damping factor \eqref{lbd2} to $\gamma\le0.8$ and use of rational coefficients finally produced a standard method with angle $\alpha=83.74$ for the triplet \texttt{AP4o33vs} for {\bf s}mooth grids.
Here, the algebraic restrictions from Lemma~\ref{Lnoderest} lead to quite long algebraic expressions for all coefficients like
\begin{align}\label{kn4o33vs}
 \cc\T=\left(\frac{144997}{389708},\frac{73}{748},\frac{77297572}{117896267},1\right)
  \doteq(0.37207,0.09759,0.65564,1),
\end{align}
with non-monotonic nodes.
The quadrature weights in $K$ are positive again in the interval $[0.11,0.43]$.
Since the coefficients of the triplet have very long rational representations, we will present them in double-precision only in Appendix~A2.
By a proper choice of $\hat b_{44}(\sigma)$ and a simplified weight matrix $W$, the interval of zero-stability in \eqref{GlmStab} becomes $[\underline\sigma,\bar\sigma]=[0.65,1.80]$.
The error constants, $err_3\doteq0.051$ and $err_3^\dagger\doteq0.032$, are larger than for \texttt{AP4o33vg}, see Table~\ref{TPT}.
Here, the form of the local errors in \eqref{lokFRestv}, \eqref{lokFResta} with $r=q=3$ is slightly different,
\begin{align}\label{LF4o33vs}
 \begin{array}{rlrl}
 \tau_n^Y=&h_n^r\beta_r(\sigma_n)y^{(r)}(t_n)+O(h_n^s\|y^{(s)}(t_n)\|_{[n]}),
  &\eins\T A\beta_r(\sigma)\equiv&0,\\[1mm]
 \tau_n^P=&h_n^q\beta_q^\dagger(\sigma_{n+1})p^{(q)}(t_n)+O(h_n^s\|p^{(s)}\|_{[n]}),
 & \eins\T A\T\beta_q^\dagger(\sigma)=&O(\sigma-1).
 \end{array}
\end{align}
Since $\beta_3^\dagger$ is a smooth function near $\sigma=1$, we will show later that global order 3 still can be preserved for {\em smooth grids} satisfying $\sigma_n=1+O(h_n)$.
Such grids are obtained naturally towards the end of adaptive grid generation strategies with successive refinements for sharper tolerances.
Smooth grids are required for enhanced accuracy of the solutions only, but not by stability issues.
In practice, the interval $[\underline\sigma,\bar\sigma]=[0.5,1.8]$ of strong zero-stability should be comfortably large enough to enable the construction of coarse grids at the start of the computation.
\par
For the boundary methods, it was again possible to choose $K_0=K_N=K$ and the block triangular form seen in \eqref{RM4o33vg} where $\mu_0=5.1$, $\mu_N=2.8$, see Table~\ref{TPB}.
The coefficients are also shown with real numbers in Appendix~\ref{SApendx}.

\subsection{A more accurate method satisfying additional order conditions}\label{S4o43}
In our previous paper \cite{LangSchmitt2022b} it was observed that in many test problems the coupling between the errors in the state variable $y$ and the multiplier $p$ is rather weak for constant stepsizes.
Methods satisfying more order conditions for the state $y$ than for $p$ indeed showed higher orders of convergence for $y$.
Accordingly, we are now trying to increase the local order for the forward method beyond $r=s-1$.
By \eqref{OBDv} additional conditions have to be solved with the elements of the last columns of $\hat A$ and $\hat B$.
Concentrating on the case $s=4$ these matrices are given by \eqref{BdAd4} and \eqref{Bsigma}.
First of all, order $4$ leads to $\hat a_{14}=1$, i.e. $e_1\T\hat A=\eins\T$.
By \eqref{precons}, this leads to $e_1\T\hat B(\sigma)=\eins\T= e_1\T\PP_s$ and shows that $e_1\T\hat B(\sigma)\PP_s^{-1}=e_1\T \QQ_{s,s}=e_1\T$ and means that also $\QQ_{2,s}$ possesses Hankel form.
Hence, the conditions \eqref{KOBstrt}, \eqref{KOBend} at the boundary are solvable even with $r=s$ for any such standard method.
Now, the remaining residual of the full order condition \eqref{OBvStd} is
\begin{align}\label{OBv4}
 \hat A-\hat K\tilde E_4-\hat B(\sigma)P_4^{-1}S^{-1}
 =&
  \begin{pmatrix}
  0\\
  (\hat a_{24}-\frac34)(1-\sigma^{-4})\\
  \hat a_{34}(1-\sigma^{-5})+(2\hat a_{24}-\frac12)\sigma^{-5}-\hat k_5(3+2\sigma^{-5})\\
  \ast
   \end{pmatrix}e_4\T.
\end{align}
The last element in the column vector is not shown since it may be chosen arbitrarily with the remaining free parameter $\hat b_{44}(\sigma)$.
Unfortunately, it seems that canceling the second entry in \eqref{OBv4} uniformly in $\sigma$ with the  choice $\hat a_{24}=\frac34$ leaves no $A(\alpha)$-stable methods at all.
Hence, we resorted to using the higher order $r_1=4$ in a restricted sense by canceling the residuals \eqref{OBv4} for $\sigma=1$ only with the choice
\begin{align}\label{OBvplus}
 \hat a_{14}=1,\ \hat a_{24}=\frac14+\frac52\hat k_{5}.
\end{align}
This leaves only 5 parameters in $\hat A$.
However, using also one degree of freedom from the nodes (see Lemma~\ref{Lnoderest}) algebra software was able to solve the equations for triangular form of the matrix $A=V\mT\hat AV^{-1}$ formally.
Extended computer searches in a three-parameter set of solutions found a very attractive standard method for the triplet \texttt{AP4o43vs} based on the node vector
\begin{align}\label{node4o43}
 \cc\T=\left(\frac1{2}(7-\sqrt{29}),\frac12,\frac1{10}(3+\sqrt{29}),1\right)
\doteq(0.0807,0.5,0.8385,1).
\end{align}
We like to mention that the order 4 in the name refers to the global order of the scheme when applied to simple initial value problems without control.
Since $c_s=1$ and by \eqref{OBvplus} this is an LSRK-method by Lemma~\ref{LLSRK}.
It possesses a stability angle of $\alpha=74.015^o$, small norm $\|\bar B(1)\|_\infty\doteq1.63$, a good damping factor $|\lambda_2|\doteq 0.52$ and small error constants $err_4\doteq 3.1\cdot 10^{-3}$, $err_3^\dagger\doteq0.076$.
It also has positive quadrature weights
\begin{align}\label{Kcolsum}
 \eins\T K\doteq(0.23926055,0.50765568,0.16243097,0.090652800)>0\T.
\end{align}
This method was designed with a moderate stepsize ratio $\hat\sigma=1.2$ in the linear program \eqref{LinProgSkal} only, but it satisfies the norm estimate \eqref{GlmStab} in a larger interval.
In order to avoid algebraic numbers in the norm $\ltnorm\cdot\rtnorm$, we constructed a simple rational approximation of its weight matrix as
\begin{align}\label{WAP4o43}
 W=\begin{pmatrix}
 1& -\frac{11}6& \frac{25}6& -\frac{7}6\\[1mm]
 1& -1& -\frac{3}2& 1\\[1mm]
 1& -\frac{1}3& -\frac{3}2& -\frac{2}5\\[1mm]
 1& 0& 0& 0
 \end{pmatrix}
\end{align}
for which $W^{-1}\bar B(\sigma)W$ also has $1+(s-1)$ block diagonal form (with a new lower block $B_{se}(\sigma)$) since it preserves the essential properties that $We_1=\eins$ is the right eigenvector of $\bar B^\infty$ and $e_1\T W^{-1}=\eins\T A=e_s\T$ its left eigenvector  which means $e_s\T W=e_1\T$.
The norms of the weight matrix and its inverse are $\|W\|_\infty\le8.2$ and $\|W^{-1}\|_\infty\le2$.
With this weight matrix, the weight $W^\dagger=(AW)\mT$ for $P$, and some fine-tuning of the parameter $\hat b_{44}(\sigma)$, the stability estimates
\eqref{GlmStab} and \eqref{NormSupKv}--\eqref{Btprest} of the standard scheme with the nodes \eqref{node4o43} hold in the interval $[\underline\sigma,\bar\sigma]=[0.47,1.79]$ with $\tilde\gamma\le 0.984<1$.
In fact, the function $\sigma\mapsto\|B_{se}(\sigma)\|_\infty$ consists of two monotone parts with minimum near $\sigma=1.5$.
Such norm estimates can even be proven rigorously by algebra software since only rational numbers and square roots are involved,
by evaluating and bounding the relative maxima at the end points exactly.
For instance, at $\sigma=145/81>1.79$ and with a simple upper bound $\sqrt{29}\le727/135$, we have that
$\|B_{se}(1.79)\|_\infty<3767/3808\doteq0.989$.
\par
The additional order conditions with $r_1=4$ for the boundary methods consume more free parameters from the matrices $(A_0,K_0)$ and $(A_N,K_N)$ and the following requirements could only be be fulfilled with full $4\times4$-block form.
These criteria concern the error constants \eqref{FeKov}, \eqref{FeKoa}, non-singularity of the step equation \eqref{defmun}, and stability of the boundary steps.
We denote this Peer triplet by \texttt{AP4o43vs} (4 stages, {\bf o}rders 4,3, {\bf v}ariable stepsize, {\bf s}mooth grid).
The essential data of the triplet are presented in Table~\ref{TPT} and its coefficients in Appendix~A3.
While the norm $\ltnorm A_N^{-1} B(1)\rtnorm$ for the end step of \texttt{AP4o43vs} equals one the norm $\ltnorm A_0\mT B(1)\T\rtnorm\cong5.2$ for the starting step  is still moderate.
See Table~\ref{TPB} for further data of the boundary methods.
\par
The LSRK property of the standard method in \texttt{AP4o43vs} leads to a further benefit for the local error since the condition for super-convergence of order $s$ is satisfied for arbitrary $\sigma$.
In \eqref{SupKvs} we get here
\begin{align}\label{SupKvRK}
 s!e_1\T\hat\beta_s(\sigma)=s!\eins\T A\beta_s(\sigma)=(1-\sigma^{-s})(\eins\T A\cc^s-1)=(1-\sigma^{-s})(c_s^s-1)\equiv 0.
\end{align}
Hence, the adjoint local error corresponds to \eqref{LF4o33vs} but the error forward \eqref{lokFRestv} with $r=q=3=s-1$ is more involved, having the representation
\begin{align}\notag
 \tau_n^Y=&\;h_n^r\beta_r(\sigma_n)y^{(r)}(t_n)+h_n^s\beta_s(\sigma_n)y^{(s)}(t_n)+O(h_n^{s+1}\|y^{(s+1)}\|_{[n]}),\\\label{LF4o43vs}
 \beta_r(\sigma)=&\;O(1-\sigma),\quad \eins\T A\beta_r(\sigma)\equiv\eins\T A\beta_s(\sigma)\equiv0.
\end{align}
In principle, different negative powers of $\sigma$ appear in $\|\beta_r(\sigma)\|$.
However, all expressions vanish for $\sigma=1$ and may be simply bounded by $|1-\sigma|$ since $\sigma_n\in[\underline\sigma,\bar\sigma]$ by assumption. 
\subsection{An A-stable triplet \texttt{AP4o33va}}
Although the stability angles of the previous Peer methods are sufficient for many problems, they do not qualify for A-stability.
Since searches of A-stable standard methods of LSRK type with $q=r=s-1=3$ were not successful, we resorted to considering Peer methods of general form.
As in \cite{LangSchmitt2022b} it seems that A-stability requires entries of different signs in the coefficient matrix $K$ and the diagonals of $A$ and nodes outside of $[0,1]$.
Such methods still qualify for stiff problems as long as the diagonals of $K^{-1}A$ are positive with $\mu_n>0$ in \eqref{defmun}.
When the orders $q=r$ coincide, one would usually look for methods with error constants $err_r\cong err_q^\dagger$ of equal magnitude.
However, for A-stable methods the size of the adjoint error constant $err_q^\dagger$ seems to be a bottleneck and we minimized the forward constant $err_r$ with higher priority.
The result was the following {\bf A}-stable method \texttt{AP4o33va} for {\bf v}ariable stepsizes based on the nodes
\begin{align*}
 \cc\T=\left(0,\frac{53}{34},\frac{6242}{30453},\frac{298}{153}\right)
  \doteq(0,1.5588,0.2050,1.9477)
\end{align*}
with a forward error $err_3\doteq 0.013$ and a larger adjoint error $err_3^\dagger\doteq0.88$.
Since the left eigenvector of $\bar B^\infty$ is not as simple for \texttt{AP4o33va} as  $\eins\T A=e_s\T$ for LSRK methods and now has long rational expressions, there is no simple rational approximation for the weight matrix $W$.
With the weight given in Appendix~A4 the method is uniformly zero-stable in the norm \eqref{GlmStab} for $\sigma\in[0.61,1.52]$, with a good damping factor $|\lambda_2|\doteq 0.29$, see also Table~\ref{TPT}.
Although $\ltnorm\bar B(1)\rtnorm=1$ holds by design, the row sum norm $\|\bar B(1)\|_\infty\doteq 20.1$ is rather large and may lead to some susceptibility for rounding errors if sharp tolerances are applied.
The size of this norm may be due to the necessity for large entries in the last row of $\hat B(\sigma)$ which seem to be required by A-stability, see Appendix~A4.
The local errors \eqref{lokFRestv}, \eqref{lokFResta} with $r=q=3=s-1$ for the standard method satisfy,
\begin{align}\label{LF4o33va}
 \begin{array}{rlr}
 \tau_n^Y=&h_n^r\beta_r(\sigma_n)y^{(r)}(t_n)+O(h_n^s\|y^{(s)}(t_n)\|_{[n]}),
  &\eins\T A\beta_r(\sigma)= O(1-\sigma),\\[1mm]
 \tau_n^P=&h_n^q\beta_q^\dagger(\sigma_{n+1})p^{(q)}(t_n)+O(h_n^s\|p^{(s)}\|_{[n]}),
 & \eins\T A\T\beta_q^\dagger(\sigma)=O(1-\sigma).
 \end{array}
\end{align}
The additional restrictions again seem to prohibit the existence of boundary methods with some reduced block structure and $K_0^{-1}A_0,\,K_N^{-1}A_N$ are full matrices.
As for the standard method the adjoint error constants are quite large for these steps, too, see Table~\ref{TPB}.
Although the row resp. column sum norms of the exceptional stability matrices are also quite large, their stability norms still are quite moderate $\ltnorm A_0\mT B(1)\T\rtnorm\doteq2.7$, $\ltnorm A_N^{-1}B(1)\rtnorm\doteq2.2$.
\begin{table}
\centerline{\begin{tabular}{|l|c|c|c|c|c|c|c|c|}\hline
  triplet &$s,(r,q)$ & $r_1$ & nodes & $\alpha$ & $\sigma\in$&$|\lambda_2|$&$err_{r_1}$&$err_q^\dagger$\\\hline
  AP4o33vg&$4,(3,3)$& 3&$[0,1]$&$61.59^o$&$[0.57,1.80]$&0.31&9.8e-3&9.8e-3\\
  AP4o33vs&$4,(3,3)$& 3&$(0,1]$&$83.74^o$&$[0.65,1.80]$&0.80&5.1e-2&3.2e-2\\
  AP4o43vs&$4,(3,3)$& 4&$(0,1]$&$74.01^o$&$[0.47,1.79]$&0.52&3.1e-3&7.6e-2\\
  AP4o33va&$4,(3,3)$& 3&$[0,1.95]$&$90^o$&$[0.61,1.52]$&0.29&1.3e-2&8.8e-1\\\hline
\end{tabular}}
\caption{Essential properties of the standard methods of Peer triplets.
}\label{TPT}
\end{table}
\begin{table}
\centerline{\begin{tabular}{|l|c|c|c|c|c|c|c|c|c|c|}\hline
 &\multicolumn{5}{c|}{Starting method}&\multicolumn{5}{c|}{End method}\\\cline{2-11}
  triplet & blksz& $q_b$ & $\mu_0$&$err_{r_1,0}$&$err_{q,0}^\dagger$& blksz&$q_b$ & $\mu_N$ & $err_{r_1,N}$&$err_{q,N}^\dagger$\\\hline
  AP4o33vg& 3+1&$\infty$ & 2.74 &1.1e-2 &8.2e-3 & 1+3&$\infty$  & 2.74&4.4e-2 &8.2e-3\\
  AP4o33vs& 3+1&$\infty$ & 5.18 &9.4e-3 &2.1e-2 & 1+3&$\infty$ & 2.84&5.4e-2 &2.7e-2\\
  AP4o43vs& 4 &2& 3.73& 1.2e-3& 8.4e-2& 4 &2& 2.93 & 3.4e-3 & 7.2e-3\\
  AP4o33va& 4 &2& 1.81& 3.1e-2& 0.77& 4 &2& 0.67& 5.6e-2& 1.17\\\hline
\end{tabular}}
\caption{Properties of the boundary methods of Peer triplets.
}\label{TPB}
\end{table}

\section{Convergence}\label{sec:convergence}
Since convergence for variable stepsizes requires more advanced arguments than in our previous papers \cite{LangSchmitt2022a,LangSchmitt2022b}, we prove it here again.
An important difficulty comes from the two-step form of the Peer methods leading to $\sigma$-dependent stability matrices in the numerical solution of the boundary problem \eqref{RWPy}, \eqref{RWPp}.
\begin{table}
\centerline{\begin{tabular}{|cl|c|cc|}
 \hline
 &&state&\multicolumn{2}{c|}{adjoint}\\\cline{3-5}
 &Time steps&order $r<s$&order $q<s$,& $q_b$\\\hline
 (a)&Start, n=0&\eqref{OBStrtv}  &\eqref{OBast},\,\eqref{A0lsb},&\eqref{OBU} \\\cline{2-5}
 (b)&Standard, $1\le n<N$&\eqref{OBvStd}&\eqref{OBast}&\\
 (c)&Super-convergence&\eqref{SupKvv}=0&\eqref{SupKva}=0&\\\cline{2-5}
 (d)&matching condition&\multicolumn{3}{c|}{\eqref{eAe}}\\\cline{2-5}
 (e)&last step, $n=N$&\eqref{OBvStd}, \eqref{EExtrp}&\eqref{OBaRB},& \eqref{OBU}
 \\\hline
\end{tabular}}
\caption{Combined order conditions for the peer triplets.}
\label{TOrd}
\end{table}
\par
Inheriting the index range from the grid, we denote the numerical solution by $Y=\big(Y_n\big)_{n=0}^N$, $P=\big(P_n\big)_{n=0}^N$, and the vector of exact solutions in boldface ${\bf y}=\big({\bf y_n}\big)_{n=0}^N$, ${\bf p}=\big({\bf p_n}\big)_{n=0}^N$, where ${\bf y_n}=\big(y^\star(t_{nj})_{j=1}^s\big)$, ${\bf p_n}=\big(p^\star(t_{nj})\big)_{j=1}^s$, $0\le n\le N$.
The global errors are indicated with checks, $\check Y=Y-{\bf y}$, $\check P=P-{\bf p}$.
For ease of writing, we also introduce the combined error vector $\check Z=\big(\check Y\T,\check P\T)\T$.
For the following analysis of the discrete boundary value problem, we multiply the equations \eqref{dBVPY0}, \eqref{dBVPYn} by appropriate inverses $A_n^{-1}$, and multiply \eqref{dBVPPn}, \eqref{dBVPPN} likewise by $A_n\mT$.
Comparing them with the residuals in these equations when the exact solution ${\bf z}$ is the argument, we obtain the following error equations:
\begin{align}\label{errY}
 \check Y_n-\bar B_n\check Y_{n-1}=&\;R_n^Y(\check Z)-\tau_n^Y,\ n=0,\ldots,N,\\\notag
  R_n^Y(\check Z):=&\;h_n\bar K_n\big(G({\bf z}_n+\check Z_n)-G({\bf z}_n)\big),
 \\\label{errP}
 \check P_n-\tilde B_{n+1}\T\tilde P_{n+1}=&\;R_n^P(\check Z)-\tau_n^P,\ n=0,\ldots,N-1,\\\notag
 R_n^P(\check Z):=&\;-h_nA_n\mT\big(\Phi({\bf y}_n+\check Y_n,K_n\T({\bf p}_n+\check P_n))-\Phi({\bf y}_n,K_n\T{\bf p}_n)\big),
\end{align}
Here, $\tau_n^Y,\tau_n^P$ are the truncation errors computed in \eqref{lokFRestv}, \eqref{lokFResta} and $\bar K_n=A_n^{-1}K_n$.
The error equation for the starting step is covered by \eqref{errY} by setting $\bar B_0=0$.
For simplicity, the adjoint boundary condition \eqref{dBVPPN} is only considered for polynomial objective functions of degree 2 at most.
Here, the error equation reads
\begin{align}\label{errRB}
 \check P_N-\left((\eins w\T)\otimes \nabla_{yy}C\right)\check Y_N=&\;R_N^P(\tilde Z)-\tau_N^P,\\\notag
 R_N^P(\check Z):=&\;-h_NA_N\mT\big(\Phi({\bf y}_N+\check Y_N,K_N\T({\bf p}_N+\check P_N))-\Phi({\bf y}_N,K_N\T{\bf p}_N)\big).
\end{align}
Combining all equations leads to a nonlinear system for the global error of the form
\begin{align}\label{GlFeGl}
 \MM_0\check Z=-\tau+R(\check Z)
 :=\begin{pmatrix}
 -\tau^Y+R^Y(\check Z)\\-\tau^P+R^P(\check Z)
 \end{pmatrix}.
\end{align}
The matrix $\MM_0$  contains all parts on the left-hand sides of \eqref{errY}--\eqref{errRB} being independent of the stepsizes $h_n$.
It has block triangular form (see also \cite{LangSchmitt2022a,LangSchmitt2022b})
\begin{align}\label{MM0}
 \MM_0=\begin{pmatrix}
  M_{11}\otimes I_m&0\\M_{21}\otimes\nabla_{yy}\CC&M_{22}\otimes I_m
 \end{pmatrix},\,
  \MM_0^{-1}=\begin{pmatrix}
  M_{11}^{-1}\otimes I_m&0\\-M_{22}^{-1}M_{21}M_{11}^{-1}\otimes\nabla_{yy}\CC&M_{22}^{-1}\otimes I_m
 \end{pmatrix}.
\end{align}
Due to the two-step form of the Peer methods, the matrix $M_{11}$ is lower block triangular and $M_{22}$ is upper block triangular with identities $I_s$ in the main diagonal and the stability matrices $\bar B_n:=A_n^{-1}B_n$, $1\le n\le N$, respectively $\tilde B_{n+1}\T:=A_n\mT B_{n+1}\T,$ $0\le n<N$, in the first off-diagonals.
Of great interest for the analysis is the precise form of their  inverses which is explicitly given by (\cite{LangSchmitt2022a})
\begin{align}\label{Invnk}
 (M_{11}^{-1})_{nk}=\bar B_n\cdots\bar B_{k+1},\ k<n,\quad
 (M_{22}^{-1})_{nk}=\tilde B_{n+1}\T\cdots \tilde B_k\T,\ k>n,
\end{align}
with identities in the main diagonal again.
The matrix $M_{21}$ has nontrivial entries of rank-one form in the very last block only, $(M_{21})_{NN}=-\eins w\T$.
Since $\eins$ is an eigenvector for all matrices $\bar B_n$ and $w\T A_N^{-1}=\eins\T$ by \eqref{OBaRB}, the southwest block in $\MM_0^{-1}$ is quite simple since $M_{22}^{-1}M_{21}M_{11}^{-1}=(\eins_N\otimes\eins_s)(\eins_N\otimes\eins_s)\T$ is a matrix of ones.
\par
For the error analysis it is of utmost importance that the products in \eqref{Invnk} are uniformly bounded.
Hence we need to use vector norms associated with the weighted matrix norms \eqref{GlmStab}, \eqref{Stbsig} from Section~\ref{SBMF} whose construction is part of the triplet design process itself.
The vector norms are defined by
\begin{align}\label{GNormen}
 \ltnorm Y\rtnorm:=\max_{n}\|W^{-1}Y_n\|_\infty,\;
 \ltnorm P\rtnorm:=\max_{n}\|(W^\dagger)^{-1}P_n\|_1,\;
 \ltnorm Z\rtnorm:=\max~\{\ltnorm Y\rtnorm,\ltnorm P\rtnorm\},
\end{align}
where $\|W^{-1}Y_n\|_\infty=\|(W\otimes I_m)^{-1}Y_n\|_\infty$ means the maximum norm and $\|\cdot\|_1$ the sum norm in the stage space, the norm in the state space $\R^m$ may be chosen by convenience.
\par
Applying the Banach Fixed Point Theorem to the error equations \eqref{errY}, etc, requires bounds on the derivatives of the right-hand sides $g,\phi$ of the reduced problem \eqref{RWPy}, \eqref{RWPp}.
Hence, we assume that constants $\Lambda$ exist such that
\begin{align}\label{AblMx}
 \|\nabla_y g\|\le\Lambda_y,\ \|\nabla_p g\|\le\Lambda_p,\quad
 \|\nabla_y\phi\|\le\Lambda_y^\dagger,\ \|\nabla_p\phi\|\le\Lambda_p^\dagger,
\end{align}
in some open tubular neighbourhood of the exact solution $(y,p)$ of \eqref{RWPy}, \eqref{RWPp}.
We also introduce $\Gamma:=\|\nabla_{yy}\CC\|$ which vanishes for linear objective functions.
\par
An important aspect of our error estimates is localized bounds containing terms like $h_n^r\|y^{(r)}\|_{[n]}$ which may show the potential for equidistribution of errors in practical computations.
Although the usual order conditions \eqref{OBvStd}, \eqref{OBaStd} are valid for any $r,q\le s\in\N$ the one-leg-condition \eqref{OBU} have been derived for $q_b\le2$ only in \cite{LangSchmitt2022b}.
The corresponding assumption may be a limitation for triplets where $K_0,K_N$ do not have diagonal form, see Table~\ref{TPB}.
\begin{theorem}\label{TKonv}
Assume that the Peer triplet satisfies the order conditions from the lines (a), (b), (d), (e) of Table~\ref{TOrd} with $2\le r=q<s$ and $q-1\le q_b\le2$, for arbitrary $\sigma\in\R$ and that there exists a norm such that $\ltnorm\bar B(\sigma)\rtnorm=\ltnorm\tilde B(\sigma)\T\rtnorm=1$ for  $\sigma\in [\underline{\sigma},\bar\sigma]$, $1\le n<N$, see \eqref{GlmStab}, \eqref{Stbsig}.
Let the solution of the problem \eqref{RWPy}, \eqref{RWPp} satisfy  $y,p\in C^{r}[0,T^\ast]$ and the upper bounds \eqref{AblMx} the conditions
\begin{align}\label{kontrkt}
 (\zeta_y\Lambda_y+\zeta_p\Lambda_p)T\le\frac12,\quad
 (\zeta_y^\dagger\Lambda_y^\dagger+\zeta_p^\dagger\Lambda_p^\dagger+\frac{s}2\Gamma)T\le\frac12,
\end{align}
with constants $\zeta,\zeta^\dagger$ which depend only on the coefficients of the Peer triplet.
Let the objective function $\CC$ be polynomial with degree not exceeding two.
\par
Then, applying the Peer triplet on a grid with stepsize ratios $\sigma_n\in[\underline\sigma,\bar\sigma]$,
$0<\underline\sigma<1<\bar\sigma$, and a sufficiently small $H=\max_{n}h_n$, its errors satisfy
\begin{align}\label{Febas}
 \ltnorm\check Z\rtnorm=\max\{\ltnorm Y-{\bf y}\rtnorm,\ltnorm P-{\bf p}\rtnorm\}
 \le \mu\big(\max_{n} h_n^{r-1}\|y^{(r)}\|_{[n]}
 + \max_{n} h_n^{q-1}\|p^{(q)}\|_{[n]}\big).
\end{align}
\end{theorem}
{\bf Proof:}
The error equation \eqref{GlFeGl} may be written in fixed point form
\begin{align}\label{FeFPF}
 \check Z=\Psi(\check Z):=-\MM_0^{-1}\tau+\MM_0^{-1}R(\check Z).
\end{align}
Considering Lipschitz differences of the right-hand sides $R_n^Y,R_n^P$ with sufficiently small $\tilde Z,\hat Z$, we get
\begin{align}\label{LipRnY}
 R_n^Y(\tilde Z)-R_n^Y(\hat Z)=&h_n\bar K_n\big(\nabla_Y G_n\cdot(\tilde Y_n-\hat Y_n)+\nabla_P G_n\cdot(\tilde P_n-\hat P_n)\big),
\end{align}
where, e.g., $\nabla_Y G_n$ is a block-diagonal matrix of integral means, $\int_0^1\nabla_y g\big((1-t)\tilde Y_{nj}+t\hat Y_{nj}\big)dt$, $j=1,\ldots,s$.
And for $R_n^P$ we get likewise
\begin{align}\label{LipRnP}
 R_n^P(\tilde Z)-R_n^P(\hat Z)=-h_n A_n\mT\big(\nabla_Y\Phi_n\cdot(\tilde Y_n-\hat Y_n)+\nabla_p\Phi_n\cdot K_n\T(\tilde P_n-\hat P_n)\big).
\end{align}
For the Lipschitz constant of $\Psi$, we have to consider the multiplication with $\MM_0^{-1}$, and for realistic bounds, the use of the weighted norms \eqref{GNormen} is required, which are different for $Y$ and $P$.
Since the handling of $M_{11}^{-1}\bar K_n\nabla_Y G_n$ is straight-forward, we consider the second term in \eqref{LipRnY} in more detail.
By assumption, in \eqref{Invnk} we have $\ltnorm\bar B_n\rtnorm=\|W\bar B_nW^{-1}\|_\infty=1$ with the possible exception for $n=0$ or $n=N$.
Since $B_N$ is a common factor in the last block row of of the lower block tringular matrix $M_{11}^{-1}$, whose norm may exceed one, we show the left-most factor in the following estimate separately,
\begin{align*}
 &\ltnorm M_{11}^{-1}\big(h_n\bar K_n\nabla_P G_n\cdot(\tilde P_n-\hat P_n)\big)_{n=0}^{N}\rtnorm\le\max_n\ltnorm B_n\rtnorm\sum_{k=0}^nh_k\|W^{-1}\bar K_k\nabla_P G_k(\tilde P_k-\hat P_k)\|_\infty\\
 &\le\ltnorm B_N\rtnorm\|W^{-1}\bar K\|_\infty\big(O(h_0)\|\tilde P_0-\hat P_0\|_\infty+\sum_{k=1}^Nh_k\|(\tilde P_k-\hat P_k)\|_\infty)\\
 &\le\ltnorm B_N\rtnorm\|W^{-1}\bar K\|_\infty\big(O(H)+T\max_{i,j}|w_{ij}^\dagger|)\ltnorm\tilde P-\hat P\rtnorm.
\end{align*}
In the last step the factor $\max_{i,j}|w_{ij}^\dagger|$ appears due to the switch to the 1-norm for $P$.
Defining constants $\zeta_p:=2\ltnorm B_N\rtnorm\|W^{-1}\bar K\|_\infty\max_{i,j}|w_{ij}^\dagger|$, and $\zeta_y$ in a similar way, we get the Lipschitz estimate
\begin{align}\label{LipRY}
 \ltnorm M_{11}^{-1}\big(R^Y(\tilde Z)-R^Y(\hat Z)\big)\rtnorm\le(\zeta_y\Lambda_y+\zeta_p\Lambda_p)T\ltnorm\tilde Z-\hat Z\rtnorm
\end{align}
for $H$ small enough.
We note that, like the given example $\zeta_p$, all $\zeta$-constants only depend on the coefficients of the Peer triplet.
\par
For the product of $M_{22}^{-1}$ and \eqref{LipRnP}, we obtain in a similar way
\begin{align}\label{LipRP}
 \ltnorm M_{22}^{-1}\big(R^P(\tilde Z)-R^P(\hat Z)\big)\rtnorm\le(\zeta_y^\dagger\Lambda_y^\dagger+\zeta_p^\dagger\Lambda_p^\dagger)T\ltnorm\tilde Z-\hat Z\rtnorm,
\end{align}
where the norm on the left-hand side is $\|(W^\dagger)^{-1}\cdot\|_1=\|W\T A\T\cdot\|_1$ and now $\ltnorm\tilde B_1\T\rtnorm\ge 1$ may be an exceptional factor in $\zeta_y^\dagger,\zeta_p^\dagger$.
Since $(W^\dagger)^{-1}\eins=e_1$, the southwest block in $\MM_0^{-1}$ also contributes
\begin{align}
 \sum_{n=0}^N\ltnorm(\eins_s\T\otimes \nabla_{yy}\CC)\big(R_n^Y(\tilde Z)-R_n^Y(\hat Z)\big)\rtnorm\le s\|\nabla_{yy}\CC\|(\zeta_y\Lambda_y+\zeta_p\Lambda_p)T\ltnorm\tilde Z-\hat Z\rtnorm.
\end{align}
This leads to the additional term in the second part of assumption \eqref{kontrkt}.
Summing up, by assumption \eqref{kontrkt}, we have shown that $\Psi$ is a contraction
\begin{align}\label{LipBedPsi}
 \ltnorm\Psi(\tilde Z)-\Psi(\hat Z)\rtnorm\le \frac12\ltnorm\tilde Z-\hat Z\rtnorm.
\end{align}
The truncation errors $\tau^Y,\tau^P$ presented in \eqref{lokFRestv}, \eqref{lokFResta} are multiplied in $\Psi(0)=-\MM_0^{-1}\tau$ by the block matrix from \eqref{MM0}.
Bounding this multiplication straightforward, one power of $h_n^r$ is lost and we get $\ltnorm\Psi(0)\rtnorm\le(\mu/2)\big(\max_{n} h_n^{r-1}\|y^{(r)}\|_{[n]} +\max_{n} h_n^{q-1}\|p^{(q)}\|_{[n]}\big)$ with some constant $\mu$.
By standard arguments we see that a closed neighborhood $\{\tilde Z:\,\ltnorm\tilde Z\rtnorm\le\varepsilon\}$ of the origin is mapped onto itself for $\ltnorm\MM_0^{-1}\tau\rtnorm=O(H^{r-1})\le \varepsilon/4$ and finally, that the fixed point $\check Z$ in \eqref{FeFPF} is bounded by
\begin{align}\label{PsiNull}
 \ltnorm\check Z\rtnorm\le2\ltnorm\Psi(0)\rtnorm=2\ltnorm\MM_0^{-1}\tau\rtnorm,
\end{align}
since $\ltnorm\check Z\rtnorm-\ltnorm\Psi(0)\rtnorm\le\ltnorm\check Z-\Psi(0)\rtnorm=\ltnorm\Psi(\check Z)-\Psi(0)\rtnorm\le\frac12\ltnorm\check Z\rtnorm$.
\qed
\par
With the first estimate from Theorem~\ref{TKonv}, we may now prove higher order error estimates.
According to \eqref{PsiNull} this is accomplished by refining the bound for $\ltnorm \MM_0^{-1}\tau\rtnorm$ 
employing super-convergence effects according to the conditions from line (d) in Table~\ref{TOrd}.
The situation is more difficult than in previous papers since these conditions may be fulfilled for constant stepsizes only, i.e. $\sigma=1$.
We remind that different Peer triplets correspond to different choices of the parameters $\hat a_{1s},\hat a_{s1}$ in Lemma~\ref{LSupKv}.
However, super-convergence is quite robust and is well observed in our numerical tests since good grid strategies lead to stepsize ratios $\sigma_n$ clustering near one for shrinking tolerances.
The following lemma takes account of this by considering {\em smooth grids} with $\sigma_n=1+O(h_n)$ for all stepsize ratios in the grid.
\par
Due to super-convergence, the error estimates will depend on several derivatives of the solutions.
In order to shorten these expressions, we introduce the following Matlab-type abbreviation for Sobolev semi-norms of a sufficiently smooth function $v(t)$,
\begin{align}\label{Sobsnm}
 \|v^{(r_1\pts r_2)}\|_{[n]}:=\max_{k=r_1,\ldots,r_2}\|v^{(k)}\|_{[n]},
  \quad r_1\le r_2.
\end{align}
\begin{lemma}\label{LSupKv}
Let the assumptions of Theorem~\ref{TKonv} be satisfied with $q=r=s-1$ and let $y,p\in C^{s}[0,T^\ast]$.
\par\noindent
a) For methods canceling both super-convergence conditions \eqref{SupKvv}, \eqref{SupKva} uniformly in $\sigma$ the errors are bounded by
\begin{align}\label{FeSupKv}
 \max\{\ltnorm Y-{\bf y}\rtnorm,\ltnorm P-{\bf p}\rtnorm\}\le \mu\big(\max_{n}h_n^{s-1}\|y^{(s-1:s)}\|_{[n]}
 +\max_{n}h_n^{s-1}\|p^{(s-1:s)}\|_{[n]}\big),
\end{align}
with some constant $\mu>0$.
\par\noindent
b) If all stepsize ratios of the grid also satisfy $\sigma_n=1+O(h_n)$, $1\le n\le N$, the error estimate \eqref{FeSupKv} still holds if the method cancels  \eqref{SupKvv}, \eqref{SupKva} for $\sigma=1$ only.
\end{lemma}
\par\noindent
{\bf Proof:}
The improved error estimates will follow by refining the estimates for $\MM_0^{-1}\tau$ in \eqref{PsiNull}.
Since the main difference between our first two methods concerns the adjoint error, we consider $\epsilon^P:=M_{22}^{-1}\tau^P$ in more detail.
By \eqref{lokFResta}, we have $\tau_n^P=h_n^q\beta_q^\dagger(\sigma_{n+1})p^{(q)}(t_n)+O(h_n^s\|p^{(s)}\|_{[n]})$.
Also, by \eqref{SupKva} holds $q!(\tilde B^\infty)\T \beta_q^\dagger(\sigma)=\hat a_{s1}(1-\sigma^q)\eins$ where $(\tilde B^\infty)\T =\eins\eins\T A\T$.
Separating the term with the exceptional matrix $\tilde B_N\T$, where $\ltnorm \tilde B_N\T\rtnorm$ may exceed one, we obtain with some positive constant $\mu_1$ that
\begin{align}\notag
 \ltnorm\epsilon^P_n\rtnorm =
 &\;\ltnorm\tau_n^P+\sum_{k=n+1}^{N-1}\tilde B_{n+1}\T\cdots \tilde B_k\T\tau_k^P+\tilde B_{n+1}\T\cdots \tilde B_N\T\tau_N^P\rtnorm\\\notag
 \stackrel{\eqref{Stbsig}}\le
 &\;\ltnorm\tau_n^P\rtnorm+\sum_{k=n+1}^{N-1}\ltnorm (\tilde B^\infty)\T\tau_k^P\rtnorm
 +\sum_{k=n+1}^{N-1}\ltnorm\big(\tilde B_{n+1}\T\cdots \tilde B_k\T\tau_k^P-(\tilde B^\infty)\T\big)\tau_k^P\rtnorm
 +\ltnorm\tilde B_N\T\tau_N^P\rtnorm\\\notag
 \stackrel{\eqref{Btprest}}\le
 &\;\mu_1 \left( h_n^q\|p^{(q)}\|_{[n]}+h_N^q\|p^{(q)}\|_{[N]} \right.\\\notag
 &\; +\left.\sum_{k=n+1}^{N-1}\big(|\hat a_{s1}||1-\sigma_{k+1}|h_k^q\|p^{(q)}(t_k)\|+h_k^s\|p^{(s)}\|_{[k]}+\tilde\gamma^{k-n}h_k^q\|p^{(q)}\|_{[k]}\big)\right)
 \\\label{epsP}
 \le &\;\mu_1\left(2+ \sum_{k=n+1}^{N-1}\Big(|\hat a_{s1}||1-\sigma_{k+1}|+h_k+\tilde\gamma^{k-n}\Big)
 \right)\max_nh_n^q\|p^{(q\pts s)}\|_{[n]}.
\end{align}
It is obvious that the constant in brackets is uniformly bounded for $\hat a_{s1}=0$ without restrictions on $\sigma_{k+1}$ since $\tilde\gamma<1$.
The corresponding estimate for $\epsilon^Y:=M_{11}^{-1}\tau^Y$ follows in the same way by separating the contribution from the exceptional matrix $\bar B_1$ and using \eqref{lokFRestv}, \eqref{NormSupKv} where \eqref{SupKvv} leads to a term 
$|\hat a_{1s}-1|\sum_{k=1}^n|1-\sigma_k|$, which vanishes for $\hat a_{1s}=1$.
\par
The contribution to $\MM_0^{-1}$ from the subdiagonal matrix block $M_{21}=(e_N\otimes\eins)(e_N\otimes w)\T$, which has entries in the very last diagonal block only, is
\begin{align*}
 -\big((M_{22}^{-1}M_{21}M_{11}^{-1})\otimes \nabla_{yy}\CC\big)\tau^Y
 =-\eins_{s(N+1)}\otimes\big(\nabla_{yy}\CC(w\T\otimes I_m)\epsilon_N^Y\big).
\end{align*}
It is obviously bounded by $O(\ltnorm\epsilon^Y\rtnorm)$.
Summarizing, we have that without further restrictions on all $\sigma_n$ it holds that $\epsilon^Y=O(h^q)$ only for methods with $\hat a_{1s}=1$ and $\epsilon^P=O(h^q)$ only for methods with $\hat a_{s1}=0$.
Since the smallest order dominates the final estimate \eqref{PsiNull}, only methods possessing both properties are covered by case a) of the assertion. The other methods require the strengthened assumption from part b), which yields $|1-\sigma_{k}|=O(h_k)$.
\qed
\par\noindent
\begin{remark}
a) We note that $\hat a_{1s}=1$ for all LSRK methods.
\par\noindent
b) The estimate \eqref{epsP} and the discussion for $\epsilon^Y$ show that the restriction on the stepsize ratios in part b) of the lemma may be weakened to
\[ \sum_{n=1}^N |1-\sigma_n|\le const\]
which would allow a uniformly bounded number of exceptions from the rule $\sigma_n=1+O(h_n)$.
\par\noindent
c) The assumption $\sigma_n=1+O(h_n)$ may seem quite restrictive in practice.
However, considering sequences of the form $h_{n+1}=h_n(1+\mu h_n),\,h_0,\mu>0,$ it is easy to prove inductively that
\[ h_n\ge h_0\prod_{j=1}^n\big(1+\mu h_0(1+\mu h_0)^{j-1}\big)\ge h_0(1+\mu h_0)^n \]
showing that more than an exponential increase of stepsizes is still possible.
\end{remark}
Using the detailed information on the triplets from \eqref{LF4o33vg}, \eqref{LF4o33vs}, \eqref{LF4o43vs}, \eqref{LF4o33va} and Table~\ref{TPT} and \ref{TPB}, we obtain
\begin{corol}\label{CFehl}
The error estimate \eqref{FeSupKv} of Lemma~\ref{LSupKv} holds with $r=q=3$ for the triplets
\\{a)} \texttt{APv4o33vg} for stepsize ratios $\sigma_n\in[\underline\sigma,\bar\sigma]=[0.5,1.8]$,
\\{b)} \texttt{APv4o33vs}, \texttt{APv4o43vs}, \texttt{APv4o33va} for stepsize ratios $\sigma_n=1+O(h_n)\in[\underline\sigma,\bar\sigma]$ with the intervals $[\underline\sigma,\bar\sigma]$ displayed in Table~\ref{TPT}.
\end{corol}
\begin{remark}
The last corollary~\ref{CFehl} does not contain a statement on improved convergence properties of the triplet \texttt{AP4o43vs}.
This is because the error structure \eqref{LF4o43vs} only suffices to show $\epsilon^Y=O(h^4)$ while $\epsilon^P$ remains with order 3.
As stated at the end of the proof of the last lemma, the smallest order dominates the general estimate \eqref{PsiNull}.
However, in our previous papers and our numerical tests, we observed that in many problems the coupling between the errors of the state variable $y$ and that of $p$ seems to be rather weak and test results show the higher order for the state indeed.
\par
In order to justify our claim we shortly sketch an estimate for $\epsilon^Y$.
With \eqref{LF4o43vs}, \eqref{Invnk}, $\bar B^\infty=\eins\eins\T A$ and neglecting for simplicity the case of an exceptional $\tilde B_1$, we get with positive constants $\mu$, $\mu_1$, and $\mu_2$,
\begin{align*}
 \ltnorm\epsilon_n^Y\rtnorm=&\;\ltnorm\sum_{k=0}^{n-1}\bar B_n\cdots\bar B_{k+1}\tau_k^Y+\tau_n^Y\rtnorm\\
  \le&\sum_{k=0}^{n-1}\ltnorm\bar B^\infty\tau_k^Y\rtnorm+\sum_{k=0}^{n-1}\ltnorm(\bar B_n\cdots\bar B_{k+1}-\bar B^\infty)\tau_k^Y\rtnorm + \ltnorm \tau_n^Y\rtnorm\\
 \stackrel{\eqref{LF4o43vs}}\le
 & \;\mu_1\Big(\sum_{k=0}^{n-1} \big(h_k^5\|y^{(5)}\|_{[k]}+\tilde\gamma^{n-k}(h_k^3|1-\sigma_k|\|y^{(3)}\|_{[k]}+h_k^4\|y^{(4)}\|_{[k]})\big)\\
  &\quad +h_n^3|1-\sigma_n|\|y^{(3)}\|_{[n]} + h_n^4\|y^{(4)}\|_{[n]})\Big)\\
 \le&\; \mu_2 \left( 2+\sum_{k=0}^{n-1}\big(h_k+\tilde\gamma^{n-k}\big) \right)
 \max_{n}h_n^4\|y^{(3:5)}\|_{[n]}
 \le \mu \max_{n}h_n^4\|y^{(3:5)}\|_{[n]},
\end{align*}
for smooth grids satisfying $|1-\sigma_k|= O(h_k)$ for all $k$.
This estimate shows that the standard method \texttt{AP4o43vs} possesses order 4 for pure initial value problems for the state variable without control.
\end{remark}

\section{Numerical tests}\label{sec:tests}
We present numerical results for all Peer triplets listed in Table~\ref{TPT}.
All calculations have been done with Matlab-Version R2021a, using the
nonlinear solver \textit{fsolve} to approximate the overall coupled scheme
(\ref{dBVPY0})--(\ref{dBVPPN}) with a tolerance $10^{-14}$.
To illustrate the rates of convergence and efficiency of variable stepsizes, 
we consider two nonlinear unconstrained optimal control problems with 
inherent boundary layers and known exact solutions. 

\subsection{Nonlinear problem}
The first problem reads:
\begin{align*}
\mbox{Minimize } y_3(0.5)  & \\
\mbox{subject to } \quad y_1'(t) =&\,  y_1^2(t) - 2y_1(t)y_2(t) + y_2^2(t) + \lambda u(t), \\
y_2'(t) =&\,\lambda y_2(t), \\
y_3'(t) =&\, 0.5 (y_1(t) - y_d(t))^2 + 0.5\,\alpha (u(t)-u_d(t))^2,\quad t\in(0,0.5], \\
y_1(0) =&\,2,\quad y_2(0)=1,\quad y_3(0)=0, 
\end{align*}
with $y_d(t)=\exp(\lambda t)+(1-t)^{-1}$ and $u_d(t)=\exp(\lambda t)$. Integrating the equation with
$y_3'(t)$ over
$[0,0.5]$ and replacing $y_3(0.5)$ in the objective function shows that the optimal control is of 
so-called tracking type minimizing the integral $\frac12\int_0^{0.5}\big((y_1-y_d)^2+\alpha(u-u_d)^2\big)dt$.
The adjoint equations are
\begin{align*}
p_1'(t) =&\, -2(y_1(t)-y_2(t))p_1(t) - (y_1(t)-y_d(t))p_3(t), \\
p_2'(t) =&\, -2(y_2(t)-y_1(t))p_1(t) - \lambda p_2(t),\\
p_3'(t) =&\, 0,\quad t\in(0,0.5],\\
p_1(0.5) =&\,0,\quad p_2(0.5)=0,\quad p_3(0.5)=1.
\end{align*}
Observe the trivial solution $p_3(t)=1$. Hence, the optimality condition 
\begin{align*}
\nabla_uf(y,u)^T\,p =&\, \lambda p_1 + \alpha (u-u_d)p_3 = 0
\end{align*}
can be used to calculate the optimal control from $u=u_d-\lambda p_1/(\alpha p_3)$, which
gives the boundary value problem (with eliminated control)
\begin{align*}
y_1'(t) =&\, y_1^2(t) - 2y_1(t)y_2(t) + y_2^2(t) + \lambda \left( u_d(t)-
\frac{\lambda}{\alpha}\frac{p_1(t)}{p_3(t)}\right),\\
y_2'(t) =&\, \lambda y_2(t),\\
y_3'(t) =&\, 0.5 (y_1(t) - y_d(t))^2 + 0.5\alpha\,\left(\frac{\lambda}{\alpha}\frac{p_1(t)}{p_3(t)}\right)^2,\\
p_1'(t) =&\, -2(y_1(t)-y_2(t))p_1(t) - y_1(t) + y_d(t),\\
p_2'(t) =&\, -2(y_2(t)-y_1(t))p_1(t) - \lambda p_2(t),\\
p_3'(t) =&\, 0,\quad t\in(0,0.5],\\
y_1(0) =&\,2,\; y_2(0)=1,\; y_3(0)=0,\; p_1(0.5)=0,\; p_2(0.5)=0,\; p_3(0.5)=1.
\end{align*}
The exact solutions are 
$y_1(t)=\exp(\lambda t)+(1-t)^{-1}$, $y_2(t)=\exp(\lambda t)$, $y_3(t)=0$,
$p_1(t)=p_2(t)\equiv 0$, $p_3(t)=1$. The components $y_1(t)$ and $y_2(t)$ exhibit the
same exponential boundary layer structure at $t=0$. We choose $\lambda=-50$ (mildly stiff regime) 
and set $\alpha=1$.

\subsubsection{Stepsize sequences with alternating and smoothly varying $\sigma$}
\begin{figure}[t!]
\centering
\includegraphics[width=7cm]{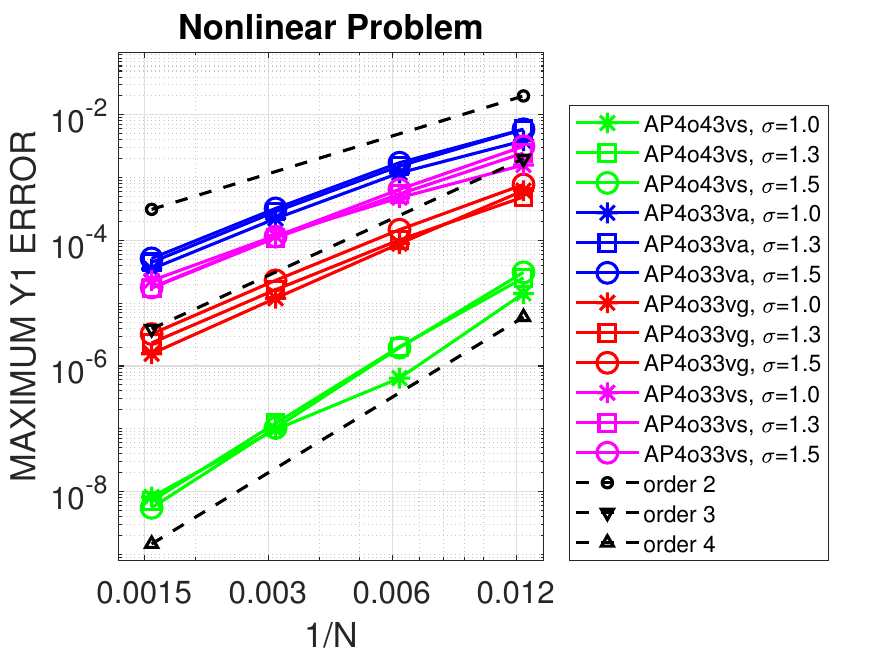}
\hspace{0.1cm}
\includegraphics[width=7cm]{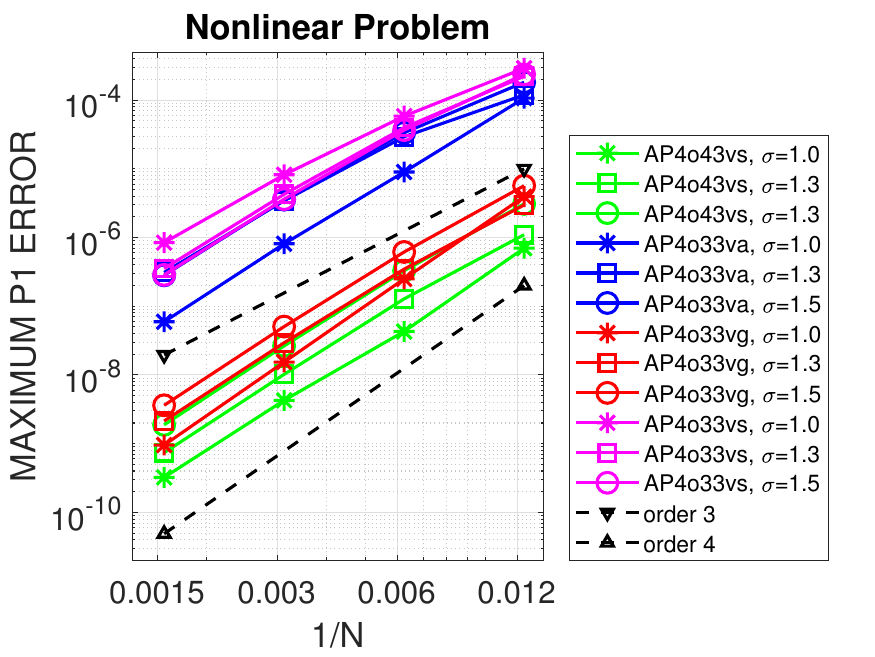}
\parbox{14cm}{
\caption{Test with alternating meshes and fixed $\sigma=1.0,\,1.3,\,1.5$. Convergence of the maximal state
errors $\|Y_{1,ni}-y_1(t_{ni})\|_\infty$ (left) and adjoint errors $\|P_{1,ni}-p_1(t_{ni})\|_\infty$ 
(right), $n=0,\ldots,N-1$, $i=1,\ldots ,4$.
}\label{fig:varsigma}
}
\end{figure}

In order to study the rates of convergence under stiffness and changing step sizes, we
first consider alternating stepsize sequences depending on a fixed parameter $\sigma$ according to
\begin{align*}
h_0=\frac{2h}{\sigma+1},\quad h_n=h_{n-1}\,\sigma^{(-1)^{n-1}},\;n=1,2,\ldots,N-1,
\end{align*}
where $h:=0.5/N$ with $N=40,80,160,320$. Note that $h_0+\ldots +h_{N-1}=0.5$ for even $N$. In 
Figure~\ref{fig:varsigma}, results for $\sigma=1.0,\,1.3,\,1.5$ are presented. Since $\sigma$ is constant, we
have $\sigma=1+\eta_n h_n$ with $\eta_n =O(h_n^{-1})$, i.e., the non-smooth case. 
Results are shown in Figure~\ref{fig:varsigma}. For the first state
component $y_1$, all methods works quite robust with respect to greater variations in the stepsize
sequence. Order four for \texttt{AP4o43vs} is clearly visible, whereas the third-order schemes reach
their order asymptotically from below. The symmetric scheme \texttt{AP4o33vg} performs significantly better than the other two third-order methods, which obviously have to pay a price for their larger stability angles $\alpha$.
This observation also applies to the first costate component $p_1$.
Here, the variations in the errors for different $\sigma$ are significant.
All methods converge with nearly fourth order for $p_1$ since $\tau^P=0$, here.
\par
\begin{figure}[t!]
\centering
\includegraphics[width=7cm]{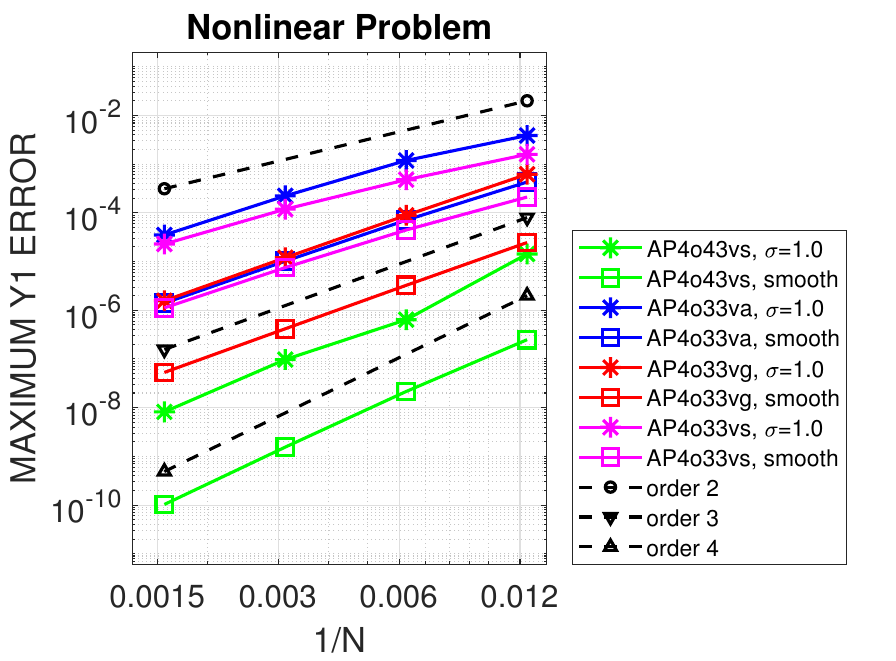}
\hspace{0.1cm}
\includegraphics[width=7cm]{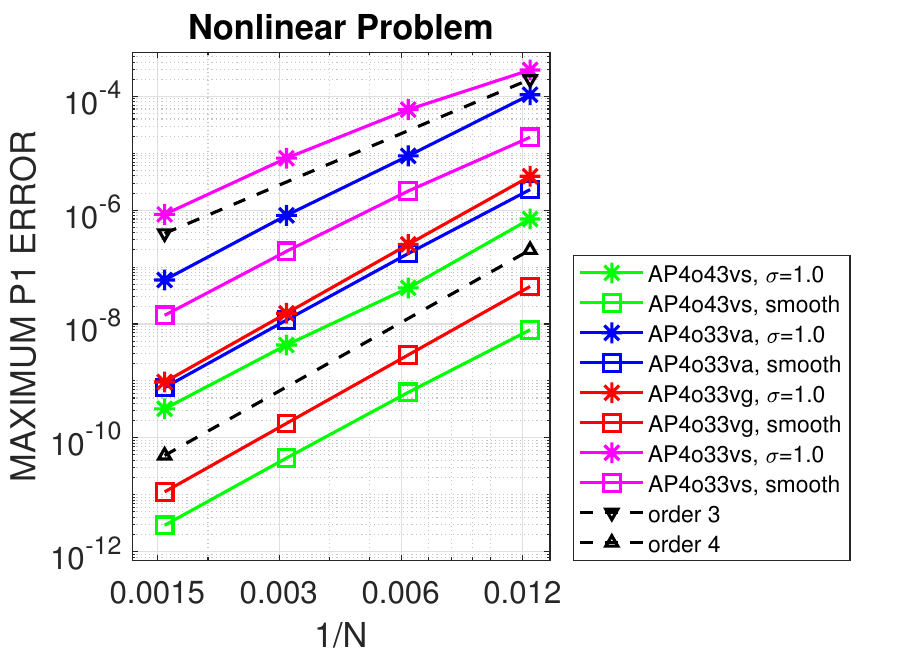}
\parbox{14cm}{
\caption{Test with smoothly varying $\sigma_n=1+0.3 h_n$ and final time $T\approx 0.22$,
results for $\sigma_n\equiv1$ are included. Convergence of the maximal state
errors $\|Y_{1,ni}-y_1(t_{ni})\|_\infty$ (left) and adjoint errors $\|P_{1,ni}-p_1(t_{ni})\|_\infty$ 
(right), $n=0,\ldots,N-1$, $i=1,\ldots ,4$.}
\label{fig:smooth}
}
\end{figure}

Next, we choose a smoothly varying stepsize sequence with $\sigma_n=1+\eta h_n$, setting 
\begin{align*}
h_0=h,\quad h_n=\frac{h_{n-1}}{1-\eta h_{n-1}},\;n=1,2,\ldots,N-1.
\end{align*}
We use $h=(0.004,0.002,0.001,0.0005)$ to mimic mesh refinement and
set $\eta=0.3$ to reach $T\approx 0.22$ for $N=40,80,160,320$.
The resulting stepsize sequence is strictly increasing which leads to grids providing 
good approximations for the stiff part $\exp(-50 t)$ of the state solution.
Results are shown in Figure~\ref{fig:smooth}. The improvements over constant time steps,
i.e. $\sigma=1$, are clearly seen. The methods \texttt{AP4o33va} and \texttt{AP4o33vs}
show their order three for the state $y_1$ for smaller time steps. Once again, 
\texttt{AP4o33vg} performs best among the third-order methods and has order three
also for larger time steps. The fourth-order \texttt{AP4o43vs} benefits from its
higher order and reaches full order four in the limit for large $N$. Order four
for the costate $p_1$ shows up for smoothly varying $\sigma_n$ by all methods, where
the symmetric \texttt{AP4o33vg} comes quite close to \texttt{AP4o43vs}.
\subsubsection{Error equidistributing stepsize sequences}
To provide an optimized stepsize sequence a priori, we select $y(t)$
for grid construction following the error equidistribution
principle introduced in \cite[Chapter 9.1.1]{AscherMattheijRussell1995}. Due to the constant solutions
$p_i(t)$, $i=1,2,3$, a well designed mesh for all components of the state vector $y$ will work fine for $p$, too. Given an asymptotic behaviour of the global error, $\|y(t_n)-Y({t_n})\|_2=O\left( h_n^r\|y^{(r)}(t_n)\|_2\right)$ with a numerical approximation $Y$ and a global convergence order $r>0$, a mesh density function 
$\psi(t):=\|y^{(r)}{(t)}\|_2^{1/r}$ is defined. Such a function has the property
\begin{align}
\int_{t_n}^{t_{n+1}} \psi(t)\,dt \approx h_n\psi(t_n) = 
h_n \| y^{(r)}(t_n)\|_2^{1/r}\approx C \|y(t_n)-Y({t_n})\|_2^{1/r}.
\end{align}
The equidistribution of the global errors over a mesh $0<t_1<\ldots <t_{N+1}=T=0.5$ requires
\begin{align}
\int_{t_n}^{t_{n+1}} \psi(t)\,dt =&\, \frac{1}{N+1} \int_0^T \psi(t)dt = const.
\end{align}
For our case, we set
\begin{align}\label{density_quad}
\psi(t):=&\,\|y^{(r)}{(t)}\|_2^{1/r}=
\left(\left(\lambda^r\exp(\lambda t)+r!\,(1-t)^{-(r+1)}\right)^2+
\left(\lambda^r\exp(\lambda t)\right)^2\right)^{\frac{1}{2r}}
\end{align}
with $r=4$ for the fourth-order method \texttt{AP4o43vs} and $r=3$ for the third-order schemes.
The Euclidean norm is used here in order to obtain a smooth function $\psi$ facilitating the 
solution of the following problem.
\par
\begin{figure}[t!]
\centering
\includegraphics[width=9cm]{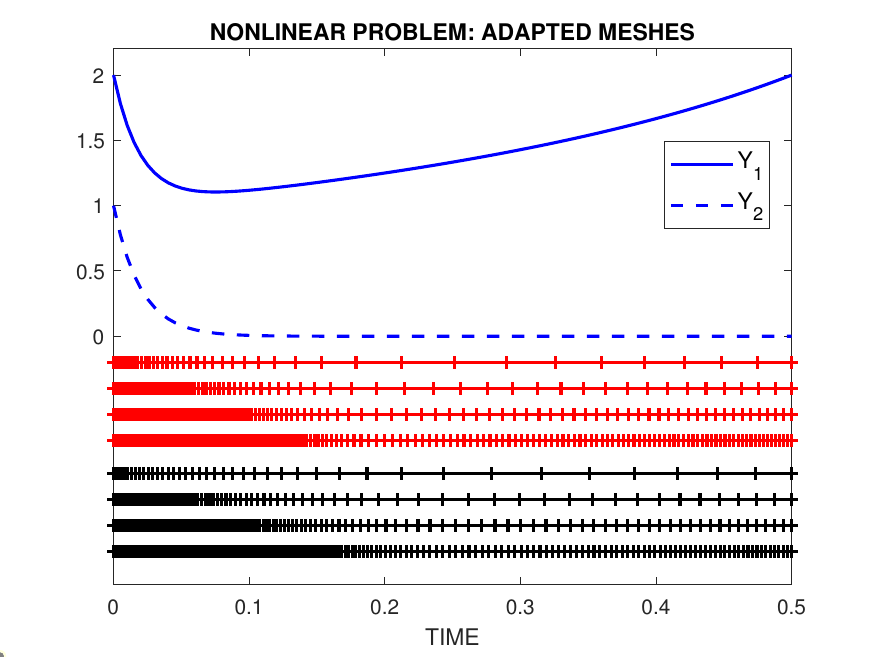}
\parbox{14cm}{
\caption{Functions $y_1(t)=(\exp(\lambda t)+(1-t)^{-1}$, $y_2(t)=\exp(\lambda t)$ 
with $\lambda=-50$ and adapted meshes for $r=3$ (red), $r=4$ (black) and $N+1=40,80,160,320$
time intervals (from top to bottom).}
\label{Fig:OptMeshesQuad50}
}
\end{figure}

 
An efficient way to derive an optimized mesh is to define a continuous node distribution $x(\xi)$
as solution of a second-order boundary value problem \cite[Chapter 2.2.2]{HuangRussell2011}
\begin{align}\label{mesh_bvp}
\partial_\xi \left( \psi(x(\xi)) \partial_\xi x(\xi) \right) =&\;0,
\quad x(0)=0, \; x(0.5)=0.5, \; \xi\in [0,0.5].
\end{align}
To numerically solve this problem, we apply uniform linear finite elements with $N$ inner nodes $\xi_n$, $n=1,\ldots,N$ and a pseudo-timestepping scheme. Finally, we set $t_n:=x(\xi_n)$.
We like to note that this procedure may also be applied with density functions $\psi$ based on error {\em estimates}.
However, this topic will be considered in a forthcoming paper.
The optimized meshes for $N+1=40,80,160,320$ are shown in Figure~\ref{Fig:OptMeshesQuad50}. The corresponding minimum and
maximum $\sigma$-values and maximum $|\eta|$-values in $\sigma_n=1+\eta_nh_n$ for this grid sequence, denoted by 
$\underline\sigma$, $\bar\sigma$, $\bar\eta$, are
\begin{align*}
\underline\sigma=(0.93,0.96,0.97,0.98),&\; \bar\sigma=(1.31,1.23,1.15,1.09),
\;\bar\eta=(16.6,16.7,18.2,19.6),\;\;\text{for }r=3,\\
\underline\sigma=(0.94,0.97,0.98,0.99),&\; \bar\sigma=(1.24,1.14,1.08,1.04),
\;\bar\eta=(12.5,12.5,12.5,12.5),\;\;\text{for }r=4.
\end{align*}
Obviously, 
for $r=3$, the maximum constants $\bar\eta$ increase only slightly in each refinement step, characterizing a smooth stepsize sequence. A closer inspection reveals that $\max_n |\eta_n|<20.4$ for arbitrary large $n$. For $r=4$, we find
$\max_n |\eta_n|\le 12.5$.

\begin{figure}[t!]
\centering
\includegraphics[width=7.3cm]{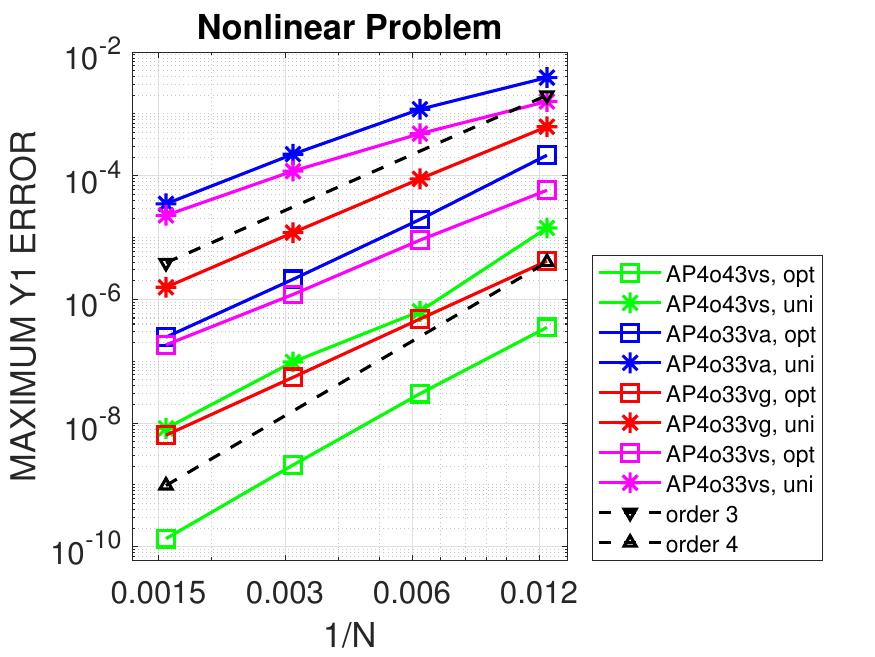}
\hspace{0.1cm}
\includegraphics[width=7.3cm]{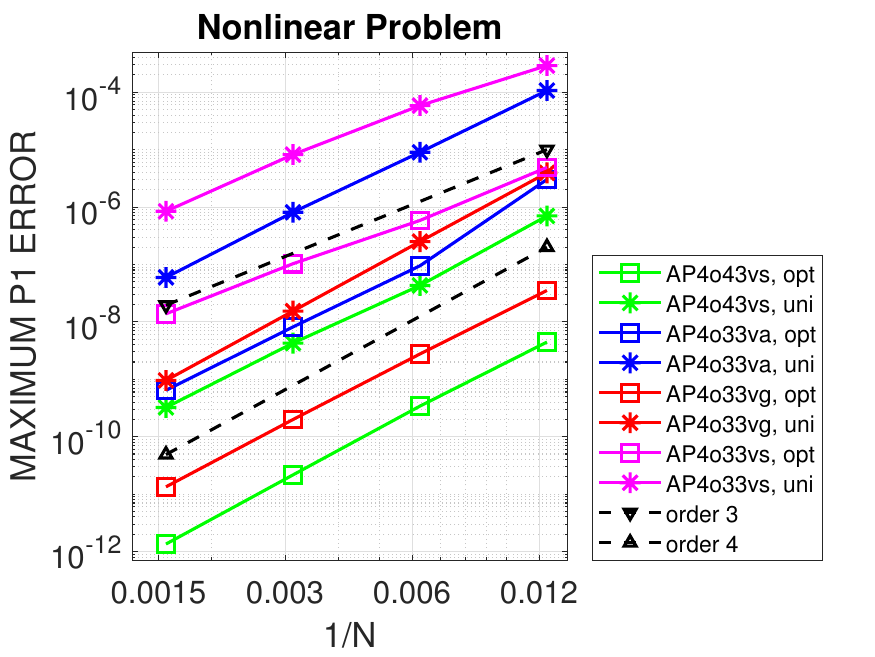}
\parbox{14cm}{
\caption{Test with adaptively chosen stepsizes, results for uniformly refined stepsizes 
are included. Convergence of the maximal state
errors $\|Y_{1,ni}-y_1(t_{ni})\|_\infty$ (left) and adjoint errors 
$\|P_{1,ni}-p_1(t_{ni})\|_\infty$ (right), $n=0,\ldots,N$, $i=1,\ldots ,4$.}
\label{fig:adapt_quad}
}
\end{figure}
Numerical results for the maximum state and adjoint errors 
of the first components $y_1$ and $p_1$ are presented in Figure~\ref{fig:adapt_quad}.
First of all, we observe for the state $y_1$ that the accuracy of all methods is improved
by two orders of magnitude versus uniform refinement. This nicely demonstrates the potential 
of global error control by adaptive variable time steps. The achieved accuracy is still one
order of magnitudes better than for the meshes in the former test with smoothly varying $\sigma$,
since the latter one could be understood as a first improvement over uniform step sizes. The third-order methods 
show their order with clear advantage for \texttt{AP4o33vg}. The fourth-order method
\texttt{AP4o43vs} asymptotically reaches order four as predicted by our convergence theory 
for smooth stepsize sequences. It beats the other triplets and the absolute values of the
errors and therefore the quality of the numerical solutions are very impressive. The
convergence orders for the costate $p_1$ are between three and four. The errors are much
smaller than for $y_1$. Once again, the symmetric \texttt{AP4o33vg} performs remarkably well
and quite close to \texttt{AP4o43vs}. 

\subsection{The Catenary problem}
The second problem is a modification of the \glqq Linea Catenaria\grqq{} of the hanging rope
\cite[Chapter 1.3, Exercise 8]{HairerNorsettWanner2006}. The well-known solution $y(t)=\cosh(a_1t+a_2)/a_1$
solves a second order differential equation, which we transform to a system of first order. The
problem formulated as an optimal control problem of tracking type reads:
\begin{align*}
\mbox{Minimize } y_3(2)  & \nonumber \\
\mbox{subject to } \quad y_1'(t) =&\,  y_2(t),\\
y_2'(t) =&\,0.5\, a_1 \left( 1 + y_2^2(t)\right)^{1/2} + u(t),\\
y_3'(t) =&\,0.5\, (y_1(t) - y_d(t))^2 + 0.5\,(u(t)-u_d(t))^2,\quad t\in(0,2],\\
y_1(0) =&\,\cosh(a_2)/a_1,\quad y_2(0)=\sinh(a_2),\quad y_3(0)=0,
\end{align*}
with $y_d(t)=\cosh(a_1t+a_2)/a_1$ and $u_d(t)=0.5\,a_1\cosh(a_1t+a_2)$. The optimality condition
\begin{align*}
\nabla_uf(y,u)\T\,p =&\, p_2 + (u-u_d)p_3 = 0
\end{align*}
with the adjoint variables $p=(p_1,p_2,p_3)\T$ can be again used to calculate the optimal 
control from $u=u_d-p_2/p_3$, which together with the adjoint equations gives the boundary 
value problem (with eliminated control)
\begin{align*}
y_1'(t) =&\, y_2(t),\\
y_2'(t) =&\, 0.5\,a_1\left( 1+y_2^2(t)\right)^{1/2} + u_d(t) - \frac{p_2(t)}{p_3(t)}\\
y_3'(t) =&\, 0.5\,(y_1(t) - y_d(t))^2 + 0.5\,\left(\frac{p_2(t)}{p_3(t)}\right)^2,\\
p_1'(t) =&\, -\left( y_1(t)-y_d(t)\right)p_3(t),\\
p_2'(t) =&\, -p_1(t) - 0.5\,a_1\left( 1+y_2^2(t)\right)^{-1/2}y_2(t)p_2(t),\\
p_3'(t) =&\, 0,\quad t\in(0,2],\\
y_1(0) =&\,\cosh(a_2)/a_1,\; y_2(0)=\sinh(a_2),\; y_3(0)=0,\; p_1(2)=0,\; p_2(2)=0,\; p_3(2)=1.
\end{align*}
The exact solutions are 
$y_1(t)=\cosh(a_1t+a_2)/a_1$, $y_2(t)=\sinh(a_1t+a_2)$, $y_3(t)=0$,
$p_1(t)=p_2(t)\equiv 0$, $p_3(t)=1$. The components $y_1(t)$ and $y_2(t)$ exhibit 
exponential boundary layer structures at $t=0$ and $t=2$. We choose $a_1=10$ and
$a_2=-10$ (mildly stiff regime).

\begin{figure}[t!]
\centering
\includegraphics[width=9cm]{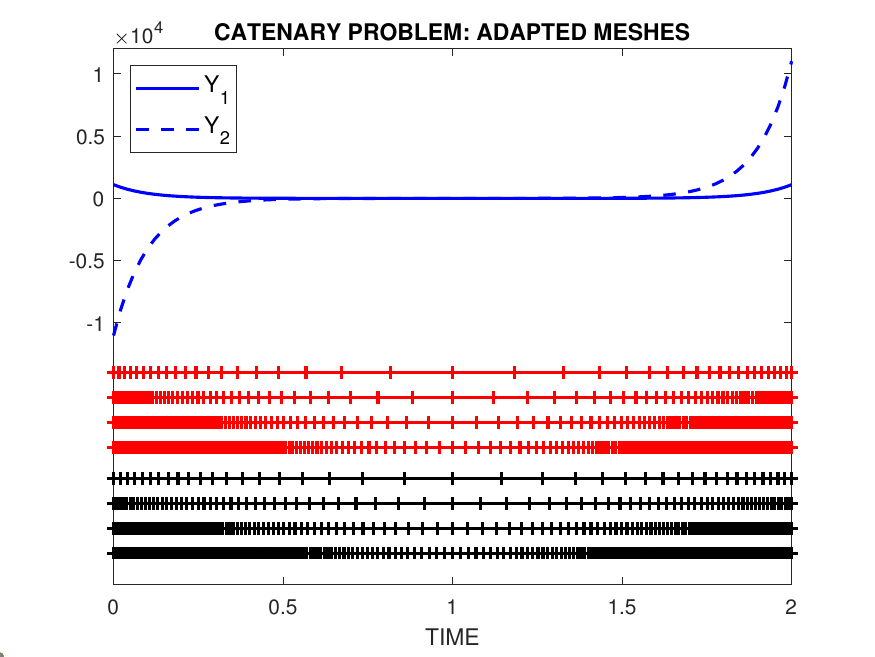}
\parbox{14cm}{
\caption{Functions $y_1(t)=0.1\cosh(10t-10)$, $y_2(t)=\sinh(10t-10)$ 
and adapted meshes for $r=3$ (red), $r=4$ (black) and $N+1=40,80,160,320$
time intervals (from top to bottom).}
\label{Fig:OptMeshesCosh10}
}
\end{figure}
\begin{figure}[t!]
\centering
\includegraphics[width=7cm]{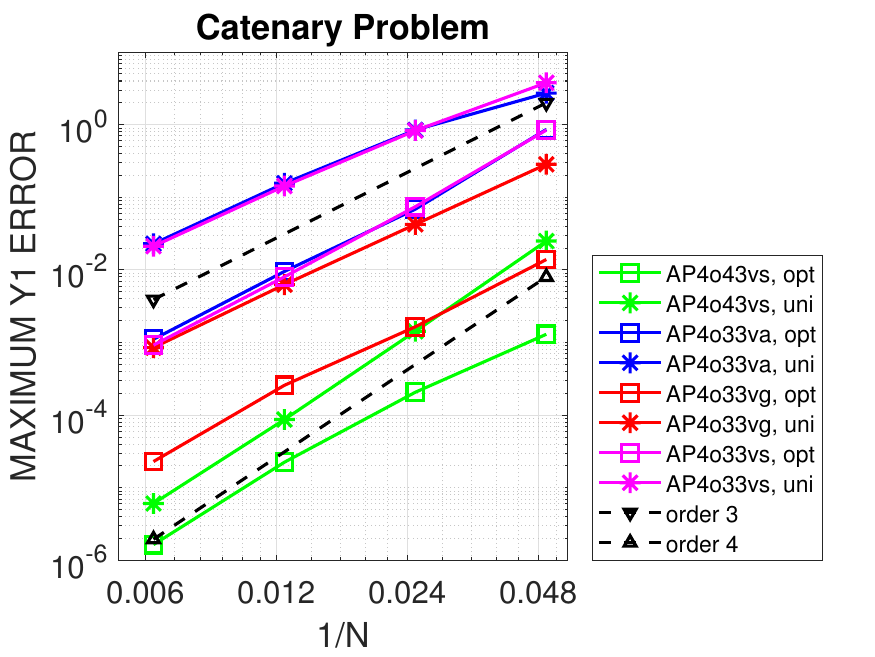}
\hspace{0.1cm}
\includegraphics[width=7cm]{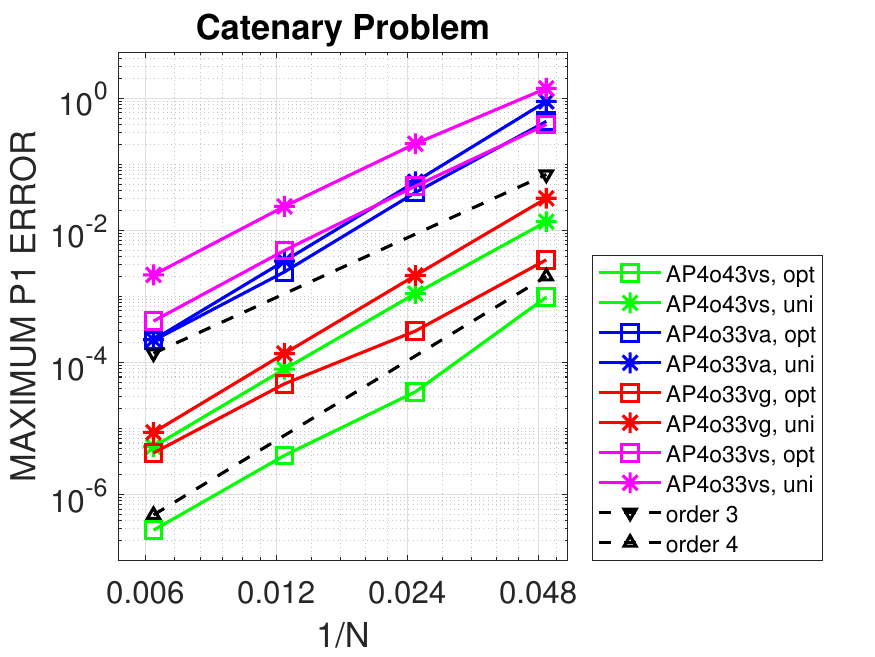}
\parbox{14cm}{
\caption{Test with adaptively chosen stepsizes, results for uniformly refined stepsizes 
are included. Convergence of the maximal state
errors $\|Y_{1,ni}-y_1(t_{ni})\|_\infty$ (left) and adjoint errors 
$\|P_{1,ni}-p_1(t_{ni})\|_\infty$ (right), $n=0,\ldots,N$, $i=1,\ldots ,4$.}
\label{fig:adapt_cosh}
}
\end{figure}
We will show results for optimized meshes adapted to $y(t)$. 
Here, the mesh density function is defined as
\begin{align}\label{density_cosh}
\psi(t):=&\,\|y^{(r)}{(t)}\|_2^{1/r}=
\left(\left( a_1^{r-1}\cosh^{(r)}(a_1t+a_2) \right)^2+
\left( a_1^r\sinh^{(r)}(a_1t+a_2)\right)\right)^{\frac{1}{2r}}
\end{align}
with $r=4$ for the fourth-order method \texttt{AP4o43vs} and $r=3$ for the third-order schemes.
Solving the boundary value problem \eqref{mesh_bvp} for $N+1=40,80,160,320$, we get the optimized 
meshes shown in Figure~\ref{Fig:OptMeshesCosh10}. The corresponding minimum and maximum $\sigma$-values 
and maximum $|\eta|$-values are
\begin{align*}
\underline\sigma=(0.74,0.80,0.86,0.92),&\; \bar\sigma=(1.35,1.25,1.16,1.09),
\;\bar\eta=(3.31,3.33,3.33,3.33),\;\;\text{for }r=3,\\
\underline\sigma=(0.81,0.87,0.92,0.96),&\; \bar\sigma=(1.24,1.15,1.08,1.04),
\;\bar\eta=(2.49,2.50,2.50,2.50),\;\;\text{for }r=4.
\end{align*}
The maximum values $\max_n |\eta_n|$
are nearly constant such that all stepsize sequences are perfectly smooth. 

Numerical results for the maximum state and adjoint errors 
of the first components $y_1$ and $p_1$ are presented in Figure~\ref{fig:adapt_cosh}.
The improvement of the accuracy for the state $y_1$ on adapted meshes is clearly visible
 but less pronounced compared to the first test problem since the boundary layer is less marked here.
All methods asymptotically reach their order three and four, respectively. 
\texttt{AP4o33va} and \texttt{AP4o33vs} perform equally well, whereas \texttt{AP4o33vg} is
again the best third-order method by two orders of magnitude. The improvement of 
\texttt{AP4o43vs} is less pronounced, but the scheme delivers still the best results. The
convergence orders for the costate $p_1$ is close to four. \texttt{AP4o33vg} and
\texttt{AP4o43vs} are the best methods with advantage for the latter one.

\section{Summary}\label{sec:summary}
We have constructed and analysed variable-stepsize implicit Peer triplets
which can be applied to control problems with varying dynamics in the constraints of
ODE type. 
Convergence results of s-stage Peer methods are proved for bounded or smoothly changing 
stepsize ratios. A notable theoretical result is that an LSRK-property 
(Last Stage is Runge-Kutta) leads to improved properties of the triplet. Balancing between
good stability properties and small error constants, we have constructed three 
third-order methods:
(1) the A($61.59^o$)-stable symmetric pulcherrima \texttt{AP4o33vg} with no smoothness
restrictions on the stepsize ratios and very small error constants, 
(2) the A($83.74^o$)-stable \texttt{AP4o33vs} for smooth grids and
with larger error constants,
(3) the A($90^o$)-stable \texttt{AP4o33va} for smooth grids and still moderate
error constants, especially suitable for dynamical systems with eigenvalues close to 
or on the imaginary axis. We could also find the A($74.015^o$)-stable \texttt{AP4o43vs}, 
which has order four for the state variables with a very small error constant, if applied
with smooth grids. All the methods show their theoretical order in our numerical tests. 
The higher-order \texttt{AP4o43vs} always gives the best results, but also our pulcherrima triplet
\texttt{AP4o33vg} performs surprisingly well and beats all other third-order methods.

In future work, we will discuss for three of the Peer triplets from this paper the integration 
of automatic grid construction in efficient gradient-based solution algorithms for the fully 
coupled optimal control problem with unknown control $u$ and further constraints, where positivity 
of the quadrature weights in $K$ is required.

\vspace{0.5cm}
\par
\noindent {\bf Acknowledgements.}
The first author is supported by the Deutsche Forschungsgemeinschaft
(DFG, German Research Foundation) within the collaborative research center
TRR154 {\em ``Mathemati\-cal modeling, simulation and optimisation using
the example of gas networks''} (Project-ID 239904186, TRR154/3-2022, TP B01).

\begin{appendix}
\section{Coefficients of all Peer triplets}\label{SApendx}
Although all data for the 4-stage Peer triplets from the text are known exactly with long rational or algebraic expressions we display most of them as real numbers in double precision due to space limitations.
The essential data are the node vector $\cc\T=(c_1,c_2,c_3,c_4)$, the coefficients $(A_0,K_0)$, $(A,K)$, $(A_N,K_N)$ of the starting, the standard and the end step and the free coefficients of the sparse matrix
\begin{align*}
 \hat B(\sigma)=\begin{pmatrix}
  1&1&1&\hat a_{14}\\
  0&0&0&\hat b_{24}(\sigma)\\
  0&0&0&\hat b_{34}(\sigma)\\
  \hat a_{41}&\hat b_{41}(\sigma)&\hat b_{43}(\sigma)&\hat b_{44}(\sigma)
 \end{pmatrix}
\end{align*}
depending on the stepsize ratio $\sigma$.
Although not required for the implementation of the triplets we also present the weight matrices $W$ in order to enable verification of uniform zero stability \eqref{GlmStab}.
All other coefficients may be computed by 
\begin{align*}
 a=A_0\eins,\ w=A_N\T\eins,\quad B(\sigma)=V_4\mT\hat B(\sigma)V_4,
\end{align*}
where $V_4=(\eins,\cc,\cc^2,\cc^3)$ is the Vandermonde matrix for the nodes from $\cc$.


\subsection*{A1: Coefficients of \texttt{AP4o33vg}}

\[ \cc\T=\left(0,\frac13,\frac23,1 \right),\quad
  K_0=K=K_N=\diag\left(\frac18,\frac38,\frac38,\frac18\right),\]
\[  A_0=
 \begin{pmatrix}
  \frac{49}{80} & \frac34 & -\frac3{16}& 0\\[1.5mm]
  -\frac{87}{80}& 0 & \frac9{16}& 0\\[1.5mm]
  \frac{87}{80} & -\frac94 & \frac{27}{16}& 0\\[1.5mm]
  -\frac{49}{80}& \frac32 & -\frac{33}{16}& 1
 \end{pmatrix},\;
A=\begin{pmatrix}
  1&0&0&0\\[1.5mm]
  -\frac94&\frac94&0&0\\[1.5mm]
  \frac94&-\frac92&\frac94&0\\[1.5mm]
  -1&\frac94&-\frac94&1
 \end{pmatrix},\;
 A_N=\begin{pmatrix}
  1& 0& 0& 0\\[1.5mm]
  -\frac{33}{16} & \frac{27}{16}& \frac9{16}& -\frac3{16}\\[1.5mm]
  \frac32 & -\frac94 & 0& \frac34\\[1.5mm]
  -\frac{49}{80}& \frac{87}{80}& -\frac{87}{80}&  \frac{49}{80}
 \end{pmatrix},
 \]
\begin{align*}
 \hat B(\sigma)=\begin{pmatrix}
  1&1&1&1\\[1mm]
  0&0&0&\frac1{36\sigma}\\[1mm]
  0&0&0&0\\[1mm]
  0&\frac{\sigma}{36}&\frac{\sigma}{18}&\frac{13}{1340}+\frac{\sigma^2}{20}
 \end{pmatrix},\quad
 W=\begin{pmatrix}
  1&-2&\frac{24}5&-\frac92\\[1mm]
  1&-\frac43&0&\frac32\\[1mm]
  1&-\frac23&-\frac85&\frac32\\[1mm]
  1 & 0 & 0 &0
 \end{pmatrix}.
\end{align*}
\subsection*{A2: Coefficients of \texttt{AP4o33vs}}
\begin{align*}
 \cc\T=&\left(\frac{144997}{389708},\frac{73}{748},\frac{77297572}{117896267},1\right),\\[1mm]
 K_n=&\diag\left(0.2089552772313791,0.2461266069992848,
  0.4259606950456414,0.1189574207236947\right)
\end{align*}
for $n=0,\ldots,N$,
\begin{align*}
 A_0=&\begin{pmatrix}
  2.773177556033415& -5.711973424498560&-0.4047906551114346& 0\\
-0.2775983738279357&  2.618694207814551& 0.1431328584722113& 0\\
 -5.101798226146757&  4.755733335146421&  2.836975327925722& 0\\
  2.606219043941277& -1.662454118462412& -2.575317531286499& 1
\end{pmatrix},\\
A=&\begin{pmatrix}
 0.7588470158140062&0&0&0\\
 0.4346633458753195& 0.5989561692950702&0&0\\
 -3.295204661275873&-0.3671669165116753& 2.473930545531403&0\\
  2.101694299586548&-0.2317892527833949&-2.473930545531403& 1
\end{pmatrix},
\end{align*}
\begin{align*}
&A_N=\\&\begin{pmatrix}
 0.7588470158140062& 0& 0& 0\\
 0.1098911012176018&  0.7137947386723661& 0.2912786335371730&-0.08134495825675107\\
 -1.064925547930965&  -1.155787455679128& 0.4736590838298028&  0.5586128875241437\\
  1.474979453185272&-0.01018461275608742& -1.911848510874736&  0.8430281717173012
\end{pmatrix},
\end{align*}
\begin{align*}
\hat a_{14}=&1,\;\hat b_{24}=0.02321239244678227/\sigma,
\;\hat b_{34}=0,\\
\hat a_{41}=&0.1010743874247749,\;
\hat b_{42}=\hat a_{41}+0.003586671392069201\,\sigma,\\
 \hat b_{43}=&\hat a_{41}+0.007173342784138403\,\sigma-0.002465255918355442\,\sigma^2,\\
 \hat b_{44}=&0.0078782707622298066+0.1683589306029579\,\sigma-0.1125\,\sigma^2+0.025\,\sigma^3,
\end{align*}
\begin{align*}
 W=&\begin{pmatrix}
 1& -\frac{49}3& -4& 27\\[1mm]
 1& -\frac{47}2& 33& -\frac{95}3\\[1mm]
 1& -9& -\frac{52}3& \frac{119}5\\[1mm]
 1& 0& 0& 0
 \end{pmatrix}.
\end{align*}
\subsection*{A3: Coefficients of \texttt{AP4o43vs}}
\begin{align*}
 \cc\T=&\left(\frac1{2}(7-\sqrt{29}),\frac12,\frac1{10}(3+\sqrt{29}),1\right),
\end{align*}
\begin{align*}
 &K_0=\\
 &\begin{pmatrix}
  0.5& 1& 0& 0\\
 -1.120097818618729& -3.509114262220923& 0.02331113741482591& -0.07507889931006730\\
  1.951080835579074&  6.817902173284554& 0.04964515498231075&  0.2324661353733601\\
 -1.097482134196919& -3.777428018384294& 0.04886693226626865& -0.04407123616946104
 \end{pmatrix},\\
 &A_0=\\
 &\begin{pmatrix}
 -2.258093793670717&  1.862197768561405&  0.8958960251093118& 0\\
  11.58487375982880& -4.941113522467058& -3.725846848775559& -0.02162218680256198\\
 -21.42711527957095&  7.401740825625927&  8.196612369685553&  0.2072923201571290\\
  12.10033531341286& -4.322825071720274& -5.366661546019306&  0.8143298666454331
 \end{pmatrix}
\end{align*}
\begin{align*}
K=&\diag(0.2392605543426944, 0.5076556795243664,0.1624309662178738,0.09065279991506543),\\
A=&\begin{pmatrix}
  2.932991332809296&  0& 0& 0\\
 -9.722226151163717&  2.605421230471736&  0& 0\\
  15.03085810481218& -5.510604377851853&  2.011734286390463& 0\\
 -8.241623286457758&  2.905183147380117& -2.011734286390463& 1
\end{pmatrix},
\end{align*}
\begin{align*}
 \hat a_{14}=&1,\;\hat b_{24}=0.006728479970272900/\sigma,\;\hat b_{34}=0,\\
 \hat a_{41}=&-0.4373259052924791,\;
 \hat b_{42}=\hat a_{41}+0.0007142621905395870\,\sigma,\\
 \hat b_{43}=&\hat a_{41}+0.001428524381079174\,\sigma+0.005699612131335000\,\sigma^2,\\
 \hat b_{44}=&\hat a_{41}+0.002142786571618761\,\sigma
 -0.01091141501818702\,\sigma^2+0.01709883639400500\,\sigma^3.
\end{align*} 
A weight matrix $W$ with short rationals is given by \eqref{WAP4o43}.
\begin{align*}
 &K_N=\\
 &\begin{pmatrix}
  0.3352224422310586&  0.6666666666666666& 0.25& 0\\
 -0.4081466631436265& -2.243551054735366& -0.9919828228000089& -0.01618666259097973\\
  0.7502573728050319&  4.650087123227831&  1.851360793436682&   0.05011862096669070\\
 -0.4323129259705010& -2.589251268789736& -0.9063392628643313&  0.03405764156058810
 \end{pmatrix},\\
 &A_N=\\
 &\begin{pmatrix}
  2.133506902525376& -1.201712432255361& 2.001196862539281& 0\\
 -6.352860439191028&  7.343234398037428&-8.042312319130696&-0.06486656040768594\\
  9.042449972383633& -12.89903567845361& 14.76669675894938& 0.6218769604713869\\
 -4.823096435717981&  6.757513712671541&-8.725581302357963& 0.4429895999362990  
 \end{pmatrix}.
\end{align*}
\subsection*{A4: Coefficients of \texttt{AP4o33va}}
\begin{align*}
 \cc\T=\left(0,\frac{53}{34},\frac{6242}{30453},\frac{298}{153}\right),
\end{align*}
\begin{align*}
 K_0=&\begin{pmatrix}
 -0.07894736842105263& -0.3541666666666667& 0.8& 0\\
  -0.5092967024450286&-0.05954441426546966& 1.5& 0\\
  -0.2793212824140483&  0.4819625026869399& 0.01024569899875302& 0\\
   0.4370032270471419&  0.2582293704863321&  -1.071186604429968&\kappa_{44}^{(0)}
 \end{pmatrix},\\
 A_0=&\begin{pmatrix}
  -2.845147129315054& -0.4034338322824405& 4.858078566685144& 0\\
  -3.334526877014251&  0.1129706979359890& 3.683717732206632& 0\\
   2.756370844715334&  0.6933008415389270&-4.233985411744457& 0\\
   2.572062980946981&  0.3827708538751709&-2.987953672161328&a_{44}^{(0)}
 \end{pmatrix},
\end{align*}
where $\kappa_{44}^{(0)}=-0.13497776057693290$, $a_{44}^{(0)}=-0.2542255953866471$.
\begin{align*}
 K=&\diag(-0.4305621262329876,0.32648079224113569, 1.239059094568785,-0.1349777605769329),\\
 A=&\begin{pmatrix}
  -6.403144243666246& 0& 0& 0\\
  -6.032436530257817& 0.4188810164603250& 0& 0\\
   7.334872792461045& 0.1741541060226739&  2.017487387419302& 0\\
   4.249467800796027& 0.1925734385846475&-0.6976301724333114&-0.2542255953866471
 \end{pmatrix},
\end{align*}
\begin{align*}
 \hat a_{14}=&1.108695652173913,\;
 \hat b_{24}=-0.4962124378026289\,/\sigma,\;
 \hat b_{34}=-0.6391248143857920\,/\sigma^2,\\
 \hat a_{41}=&4.607142857142857,\;
 \hat b_{42}=\hat a_{41}-0.2679484769093443\,\sigma,\;
 \hat b_{43}=\hat a_{41}-0.5358969538186886\,\sigma,\\
 \hat b_{44}=&-\frac{2198}{55}+\frac{1607}{22}\,\sigma-\frac{147}5\,\sigma^2,
\end{align*}
\begin{align*}
 W=&\begin{pmatrix}
  1& -2& 6& 2.509523385281405\\
  1&  1.117647058823529& 12.32698961937716&  7.052310433008046\\
  1& -1.590056808852987& 1.836922096090924& 0.5167228373603610\\
  1&  1.895424836601307& 27.53940792003076&  16.07292401233786
 \end{pmatrix}
\end{align*}
\begin{align*}
 &K_N=\\&\begin{pmatrix}
 -0.4305621262329876& 0& 0& 0\\
 -0.7584777455167840& -0.4907990379996561& 0.66666666666666667& 0.9090909090909091\\
 -0.4295489737333543& 0.09171637742180421&  1.7797533837522101&-0.2028616928718752\\
  0.6522412050328370&  0.3192953012669550&-0.51790142522729154&-0.5886128416494334
  \end{pmatrix},\\
  &A_N=\\
  &\begin{pmatrix}
  -6.4031442436662458& 0& 0& 0\\
 -0.95260517222681956&   1.865037832767615& -6& -0.5259881743382867\\
   6.79592553738488869& 0.3023172450424132& 2.591223325518416& -0.1629518220426969\\
  0.525998613319805586&-0.7619282537426990& 3.676768004714683& 0.04934710726892677
  \end{pmatrix}
\end{align*}
\end{appendix}
\bibliographystyle{plain}
\bibliography{bibpeeropt}

\end{document}